\documentclass{amsart}
\usepackage{amssymb}
\usepackage[utf8]{inputenc}
\usepackage{hyperref,color,fullpage}
\usepackage{scrextend}
\usepackage{multicol}
\usepackage{subfig}
\usepackage{tikz}
\usepackage{tkz-tab}
\usetikzlibrary{automata, positioning}
\usetikzlibrary{shapes.geometric,positioning}
\tikzstyle{vertex}=[circle, draw, inner sep=0pt, minimum size=2pt]
\newcommand{\vertex}{\node[vertex]}
\usepackage{caption}
\usepackage{systeme}
\usepackage{verbatim}
\usepackage{cite}
\usepackage{bbm}
\usepackage{pgf-pie}
\usepackage{listings}
\definecolor{lightgray}{gray}{0.7}
\definecolor{midgray}{gray}{.9}
\definecolor{dkgreen}{rgb}{0,0.6,0}
\definecolor{gray}{rgb}{0.5,0.5,0.5}
\definecolor{mauve}{rgb}{0.58,0,0.82}
\definecolor{lightgray}{gray}{0.7}
\definecolor{midgray}{gray}{.9}

\definecolor{lightgray}{gray}{9}

\title{Markov models for the tipsy cop and robber game on graphs}
\author{Viktoriya Bardenova, Vincent Ciarcia, Erik Insko}
\address{Florida Gulf Coast University, 10501 FGCU Boulevard, Fort Myers, FL 33965}
\email{vbardenova@fgcu.edu, vciarcia@fgcu.edu, einsko@fgcu.edu}
\date{May 2020}

\newcommand{\ds}{\displaystyle}

\newcommand{\N}{\mathcal{N}}

\newcommand{\n}{\mathbf{n}}
\newcommand{\R}{\mathfrak{n}}
\newcommand{\G}{\mathfrak{G}}

\newcommand{\M}{\mathbf{M}}

\newcommand{\tc}{t_c}
\newcommand{\tr}{t_r}
\begin{document}

\maketitle

\begin{abstract}
In this paper we analyze and model three open problems posed by Harris, Insko, Prieto-Langarica, Stoisavljevic, and Sullivan in 2020 concerning the tipsy cop and robber game on graphs.   The three different scenarios we model account for different biological scenarios.  The first scenario is when the cop and robber have a consistent tipsiness level though the duration of the game; the second is when the cop and robber sober up as a function of time; the third is when the cop and robber sober up as a function of the distance between them.  Using Markov chains to model each scenario we calculate the probability of a game persisting through $\M$ rounds of the game and the expected game length given different starting positions and tipsiness levels for the cop and robber.   
\end{abstract}

\section{Introduction} 
The game of cops and robbers on graphs was introduced independently by Quilliot \cite{Q86} and Nowakowski and Winkler \cite{NW83}.  In this game, a cop and a robber alternate turns moving from vertex to adjacent vertex on a connected graph $G$ with the cop trying to catch the robber and the robber trying to evade the cop.  In 2014, Komarov and Winkler studied a variation of the cop and robber game in which the robber is too inebriated to employ an escape strategy, and at each step he moves to a neighboring vertex chosen uniformly at random \cite{KW14}.

In 2020, Harris, Insko, Prieto-Langarica, Stoisavljevic, and Sullivan introduced another variant of the cop and robber game that they call the tipsy cop and drunken robber game on graphs \cite{HIPSS20}. 
Each round of this game consists of independent moves where the robber begins by moving uniformly and randomly on the graph to an adjacent vertex from where he began, this is then followed by the cop moving to an adjacent vertex from where she began; since the cop is only tipsy, some percentage of her moves are random and some are intentionally directed toward the robber. 

In this paper, we generalize the work of Harris et al. to analyze the tipsy cop and tipsy robber game.
One inspiration for this study to model is the biological scenario illustrated in the YouTube video \cite{I13} \url{ https://www.youtube.com/watch?v=Z_mXDvZQ6dU} where a neutrophil chases a bacteria cell moving in random directions.  While the bacteria's movement seems mostly random, the neutrophil's movement appears slightly more purposeful but a little slower. 

Harris et al. modeled the game by having the players alternate turns. 
In this paper we use a slightly different model of the game that was suggested to us by Dr. Florian Lehner.  We let the cop and robber be any amount of tipsy and rather than assume that the players alternate turns and flip a coin to determine if the player moves randomly or not, we use a spinner wheel to determine the probability of whether the next move will be a sober cop move, a sober robber move, or a tipsy move by either player.  
We model this scenario on vertex-transitive graphs (both finite and infinite examples) and on non-vertex transitive graphs (friendship graphs) using the theory of Markov chains.  Given a set of initial conditions on each player's tipsiness and their initial distance, the questions we consider  are:

\begin{itemize}
\item What is the probability $P(i,j,\M)$ that the game, beginning in state $i$, will be in state $j$ after exactly $\M$ rounds?
\item What is the probability $\G_{\M}(d)$ that the game lasts at least $\M$ rounds if they start distance $d$ away?
\item What is the expected number $E(d)$ of rounds the game should last if they start distance $d$ away?
\end{itemize}

For some of these graphs we are able to identify when the expected capture time is finite and when the robber is expected to escape.

In Section \ref{sec:Background} we introduce our model of the tipsy cop and robber game on both vertex-transitive and non-vertex transitive graphs.
In Section \ref{sec:Gambler's ruin} we analyze the game on regular trees and compare it to the famous Gambler's ruin problem \cite{S09}.
In Section \ref{sec:General_Method} we present our general method for analyzing the game using Markov chains, and we then apply these methods to analyze the game on specific families of graphs in Sections \ref{sec:Cycle Graphs}- \ref{sec:Toroidal}.
In Section \ref{sec:Cycle Graphs} we analyze the game on cycle graphs.
In Section \ref{sec:Petersen} we analyze the game on Petersen graphs.
In Section \ref{sec:Friendship Graphs} we analyze the game on friendship graphs.
In Section \ref{sec:Toroidal} we analyze the game on toroidal grids.
In Section \ref{sec:Sober} we present a model for the game where players start off drunk and sober up as the game progresses, answering Question 6.1 from Harris et al. \cite{HIPSS20} 
Finally, in Section \ref{sec:Distance} we present a model for the game where the players' tipsiness increases as a function of the distance between them, thus providing an answer to Question 6.6 from Harris et al. \cite{HIPSS20}. While we present our methods with specific examples from each family of graphs we analyze, we also include our Cocalc (Sage) code in the Appendix so that any reader may adapt our methods to their own models.

\section{Background and introduction of our model}\label{sec:Background}
\begin{figure}[h]
\begin{tabular}{cc}

\begin{tikzpicture}[scale = 1.65]
\pie[
  color = {yellow, green, cyan, red},radius = .75,
  /tikz/nodes={text opacity=1,overlay}
] {30/r,30/c,15/$t_c$,25/$t_r$}
\draw [ultra thick][->] (0,0) -- (.0,.4743);
\end{tikzpicture}
\end{tabular}

\captionof{figure}{Probability of each move based on spinner.} \label{fig:spinner4}
\end{figure}
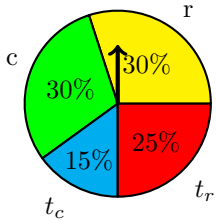

In this paper we model the tipsy cop and tipsy robber game on a graph $G$ by first placing the cop on one vertex and the robber on another vertex on the graph $G$. Rather than require the players alternate turns as in previous models of the cops and robbers game, we allow for four possible outcomes in each round of the game: a sober robber move $r$, a sober cop move $c$, a tipsy robber move $\tr$, or a tipsy cop move $\tc$, where $c+r+\tc+\tr=1$. The outcome of each round is assigned at random, perhaps by a probability spinner as depicted in Figure \ref{fig:spinner4}.

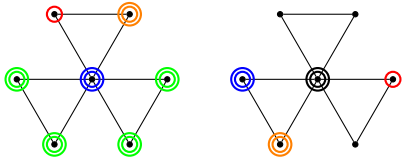
\begin{figure}[h]
\begin{tikzpicture}
\vertex[fill](a1) at (.5000000000, .8660254040) {};
\vertex[fill](a2) at (-.5000000000, .8660254040) {};
\vertex[fill](a3) at (-1., 0.) {};
\vertex[fill](a4) at (-.5000000000, -.8660254040) {};
\vertex[fill](a5) at (.5000000000, -.8660254040) {};
\vertex[fill](a6) at (1., 0.) {};
\vertex[fill](a0) at (0., 0.) {};
\draw (a1)--(a2) ;
\draw (a0)--(a2) ;
\draw (a3)--(a4) ;
\draw (a0)--(a4) ;
\draw (a3)--(a0) ;
\draw (a0)--(a1) ;
\draw (a5)--(a6);
\draw (a5)--(a0);
\draw (a0)--(a6);
\draw[orange,thick] (.5000000000, .8660254040)circle (.1cm);
\draw[orange,thick] (.5000000000, .8660254040)circle (.15cm);
\draw[red,thick] (-.5000000000, .8660254040) circle (.1cm);
\draw[green,thick] (-1., 0.) circle (.1cm);
\draw[green,thick] (-1., 0.) circle (.15cm);
\draw[green,thick] (-.5000000000, -.8660254040) circle (.1cm);
\draw[green,thick] (-.5000000000, -.8660254040) circle (.15cm);
\draw[green,thick] (.5000000000, -.8660254040) circle (.1cm);
\draw[green,thick] (.5000000000, -.8660254040) circle (.15cm);
\draw[green,thick] (1,0) circle (.1cm);
\draw[green,thick] (1,0) circle (.15cm);
\draw[blue,thick] (0,0) circle (.1cm);
\draw[blue,thick] (0,0) circle (.15cm);
\vertex[fill](a1) at (3.5000000000, .8660254040) {};
\vertex[fill](a2) at (3-.5000000000, .8660254040) {};
\vertex[fill](a3) at (3-1., 0.) {};
\vertex[fill](a4) at (3-.5000000000, -.8660254040) {};
\vertex[fill](a5) at (3.5000000000, -.8660254040) {};
\vertex[fill](a6) at (3+1., 0.) {};
\vertex[fill](a0) at (3+0., 0.) {};
\draw (a1)--(a2) ;
\draw (a0)--(a2) ;
\draw (a3)--(a4) ;
\draw (a0)--(a4) ;
\draw (a3)--(a0) ;
\draw (a0)--(a1) ;
\draw (a5)--(a6);
\draw (a5)--(a0);
\draw (a0)--(a6);
\draw[blue,thick] (3-1,0) circle (.1cm);
\draw[blue,thick] (3-1,0) circle (.15cm);
\draw[red,thick] (3+1,0) circle (.1cm);
\draw[black,thick] (3.,0)circle (.1cm);
\draw[black,thick] (3.,0)circle (.15cm);
\draw[orange,thick] (3-.5000000000, -.8660254040) circle (.1cm);
\draw[orange,thick] (3-.5000000000, -.8660254040) circle (.15cm);
\end{tikzpicture}
\captionof{figure}{Possible cop (double circles) moves on a friendship graph with 3 triangles.} \label{fig:cop3frienship}
\end{figure}

\begin{itemize}
    \item $r$ is the probability of a sober robber move (in yellow)
    \item $c$ is the probability of a sober cop (in green)
    \item $\tr$ is the probability of a tipsy move by the robber (in red)
    \item $\tc$ is the probability of a tipsy move by the cop (in blue)
\end{itemize}
If the game is played on a non-vertex transitive graph the probability of transition from one state to another depends on the starting location of the players. One such graph is a friendship graph, where $n$ copies of the cycle graph $C_3$ are joined by a common vertex (see Figure \ref{fig:cop3frienship} for an example of a friendship graph with $3$ copies of the $C_3$ cycle graph).

A vertex-transitive graph is a graph, where every vertex has the same local environment, so that no vertex can be distinguished from any other based on the vertices and edges surrounding it. On a non-vertex-transitive graph a tipsy move by the cop is not necessarily equivalent to a tipsy move by the robber, as it may increase, keep the same, or decrease the distance between the players based on each player's starting state. For example, the first part of Figure \ref{fig:cop3frienship} gives us the possible movements of the cop (in blue) if she is in the center of a friendship graph with $3$ triangles. If the cop makes a tipsy move to a green vertex (which accounts for $\frac{4}{6}$ of $\tc$ moves), the distance between her and the robber (in red) increases to $2$. If she moves tipsy to the vertex in orange the distance will not change ($\frac{1}{6}$ of $\tc$ moves), and if she moves to the vertex in red (that move can be either $\frac{1}{6}$ of $\tc$ moves or the only sober move) she will catch the robber (distance will be $0$). On the other hand, if the cop's starting position is on an outer vertex (second part of Figure \ref{fig:cop3frienship}) there are only two possible outcomes, she moves closer to the robber (in black-either $\frac{1}{2}$ of $\tc$ or the only sober move) or keeps the same distance (in orange, $\frac{1}{2}$ of $\tc$ moves). 

\hspace{1.5cm}

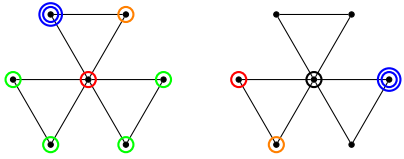
\begin{figure}[h]
\begin{tikzpicture}
\vertex[fill](a1) at (.5000000000, .8660254040) {};
\vertex[fill](a2) at (-.5000000000, .8660254040) {};
\vertex[fill](a3) at (-1., 0.) {};
\vertex[fill](a4) at (-.5000000000, -.8660254040) {};
\vertex[fill](a5) at (.5000000000, -.8660254040) {};
\vertex[fill](a6) at (1., 0.) {};
\vertex[fill](a0) at (0., 0.) {};
\draw (a1)--(a2) ;
\draw (a0)--(a2) ;
\draw (a3)--(a4) ;
\draw (a0)--(a4) ;
\draw (a3)--(a0) ;
\draw (a0)--(a1) ;
\draw (a5)--(a6);
\draw (a5)--(a0);
\draw (a0)--(a6);
\draw[orange,thick] (.5000000000, .8660254040)circle (.1cm);
\draw[blue,thick] (-.5000000000, .8660254040) circle (.1cm);
\draw[blue,thick] (-.5000000000, .8660254040) circle (.15cm);
\draw[green,thick] (-1., 0.) circle (.1cm);
\draw[green,thick] (-.5000000000, -.8660254040) circle (.1cm);
\draw[green,thick] (.5000000000, -.8660254040) circle (.1cm);
\draw[green,thick] (1,0) circle (.1cm);
\draw[red,thick] (0,0) circle (.1cm);
\vertex[fill](a1) at (3.5000000000, .8660254040) {};
\vertex[fill](a2) at (3-.5000000000, .8660254040) {};
\vertex[fill](a3) at (3-1., 0.) {};
\vertex[fill](a4) at (3-.5000000000, -.8660254040) {};
\vertex[fill](a5) at (3.5000000000, -.8660254040) {};
\vertex[fill](a6) at (3+1., 0.) {};
\vertex[fill](a0) at (3+0., 0.) {};
\draw (a1)--(a2) ;
\draw (a0)--(a2) ;
\draw (a3)--(a4) ;
\draw (a0)--(a4) ;
\draw (a3)--(a0) ;
\draw (a0)--(a1) ;
\draw (a5)--(a6);
\draw (a5)--(a0);
\draw (a0)--(a6);
\draw[red,thick] (3-1,0) circle (.1cm);
\draw[blue,thick] (3+1,0) circle (.1cm);
\draw[blue,thick] (3+1,0) circle (.15cm);
\draw[black,thick] (3.,0)circle (.1cm);
\draw[orange,thick] (3-.5000000000, -.8660254040) circle (.1cm);
\end{tikzpicture}
\captionof{figure}{Possible robber (single circles) moves on a friendship graph with 3 triangles.} \label{fig:robber3frienship}
\end{figure}

Similarly, the starting state of the robber determines the probability of moving to the next state. The first part of Figure \ref{fig:robber3frienship} depicts the possible movements of the robber (in red) if he is in the center of a friendship graph with $3$ triangles. If the robber makes a sober or tipsy (which accounts for $\frac{4}{6}$ of $\tr$ moves) move to a green vertex, the distance between him and the cop (in blue) increases to $2$. If he moves tipsy ($\frac{1}{6}$ of $\tr$ moves) to the vertex in orange the distance will not change, and if he moves to the vertex in blue (a different tipsy move $\frac{1}{6}$ of $\tr$ moves) he will get caught (distance will be $0$). On the other hand, if the robber's starting position is on an outer vertex (second part of Figure \ref{fig:robber3frienship}) there are only two possible outcomes, he moves closer to the cop (in black, $\frac{1}{6}$ of $\tr$ moves) or keeps the same distance (in orange, $\frac{1}{6}$ of $\tr$ moves or a sober move).

We will analyze the friendship graph game in more detail in Section \ref{sec:Friendship Graphs}.

The game is simpler to model on vertex-transitive graphs.  A sober cop move will always decrease the distance between the two (as the cop is chasing the robber). A sober robber move will always increase the distance between the two, or if they are at a maximum distance from each other on a finite graph, he can decide to stay in the same place. Finally, a tipsy move by either player is a random move, and therefore may increase or decrease the distance between the two.   
When modeling the game on a vertex-transitive graph, a tipsy move by the cop is equivalent to a tipsy move by the robber, so we regroup all tipsy moves together. Every move in the game is guaranteed to fall in one of those three categories hence,  $ c + r + t = 1$. 

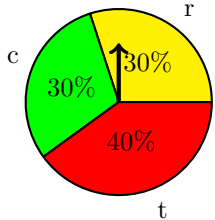
\begin{figure}[h]
\begin{tabular}{cc}
\begin{tikzpicture}[scale = 1.65]
\pie[
  color = {yellow, green, red},radius = .75,
  /tikz/nodes={text opacity=1,overlay}
] {30/r,30/c,40/t}
\draw [ultra thick][->] (0,0) -- (.0,.4743);
\end{tikzpicture}
\end{tabular}
\captionof{figure}{Probability of each move based on spinner where tipsy moves are grouped together.} \label{fig:spinner3}
\end{figure}

\begin{itemize}
    \item $r$ is the probability of a sober robber move (in yellow)
    \item $c$ is the probability of a sober cop (in green)
    \item $t$ is the probability of a tipsy move by either player (in red)
\end{itemize}
For example, if the game is played on an infinite path (Figure \ref{fig:infpathmove}), the probability that the distance between the cop and robber increases by one in a round is equal to the probability of a sober robber move, since a sober robber will always run from the cop, plus one half of the probability of a tipsy move by either player, since half of all random moves will increase the distance between them.  We will employ Markov chains and transformation matrices to model how the players move from one state to another on these graphs. 

\begin{center}
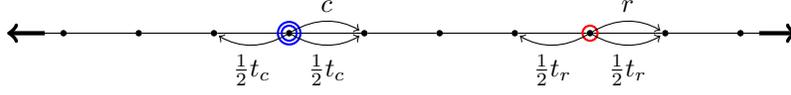
\begin{figure}[h] 
    \begin{tikzpicture} 
        \vertex[fill] (a1) at (0,0) {};
        \vertex[fill] (a2) at (1,0) {};
        \vertex[fill] (a3) at (2,0) {};
        \vertex[fill] (a4) at (3,0) {};
        \vertex[fill] (a5) at (4,0) {};
        \vertex[fill] (a6) at (5,0) {};
        \vertex[fill] (a7) at (6,0) {};
        \vertex[fill] (a8) at (7,0) {};
        \vertex[fill] (a9) at (8,0) {};
        \vertex[fill] (a10) at (9,0) {};
        
        \draw[][-] (0,0) -- (-.25,0);
        \draw [ultra thick][->] (-.25,0) -- (-.75,0);
        \draw[][-] (9,0) -- (9.25,0);
        \draw [ultra thick][->] (9.25,0) -- (9.75,0);
        \draw (a1)--(a2) (a2)--(a3) (a3)--(a4) (a4)--(a5) (a5)--(a6) (a6)--(a7) (a7)--(a8) (a8)--(a9) (a9)--(a10);
        \draw[blue,thick] (3,0) circle (.1cm);
        \draw[blue,thick] (3,0) circle (.15cm);
        \draw[red,thick] (7, 0) circle (.1cm);
        
        \draw[every loop]
            (a4) edge[bend left, auto=left] node{$\frac{1}{2}\tc$} (a3)
            (a4)edge[bend right, auto=right]node{$\frac{1}{2}\tc$} (a5)
            (a4) edge[bend left, auto=left] node {$c$} (a5)
            (a8) edge[bend left, auto=left] node{$\frac{1}{2}\tr$} (a7)
            (a8)edge[bend right, auto=right]node{$\frac{1}{2}\tr$} (a9)
            (a8) edge[bend left, auto=left] node {$r$} (a9);
    \end{tikzpicture}
\captionof{figure}{Move of a tipsy cop (blue) and tipsy robber (red)} \label{fig:infpathmove}
\end{figure}
\end{center}

\section{Regular trees and the Gambler's ruin problem}\label{sec:Gambler's ruin}
On an infinite regular tree of degree $\Delta$, the distance between the two players decreases with probability $\frac{1}{\Delta}$ and increases with probability $\frac{\Delta-1}{\Delta}$ when either player makes a tipsy move. When the cop makes a sober move, the distance always decreases, and when the robber makes a sober move the distance always increases.     
Hence, if we assume that the cop calls off the hunt when the distance between the players reaches a specified distance $d=\R$, then the Markov chain in Figure \ref{fig:Markov_trees} models the game where   the probability of the distance increasing is  $p= t \left ( \frac{\Delta-1}{\Delta} \right)+r$ and the probability of the distance decreasing is $1-p = c+ \frac{t}{\Delta}$.
\begin{center}
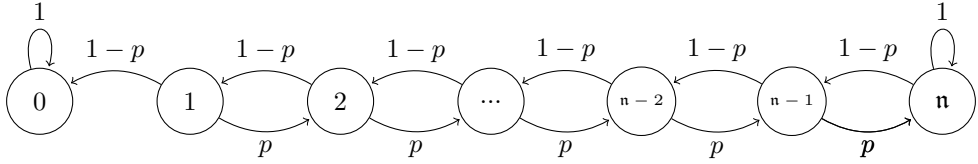
\begin{figure}[h] 
    \begin{tikzpicture} 
        \node[state] at (0,0) (0) {$0$};
        \node[state] at (2,0) (1) {$1$};
        \node[state] at (4,0) (2) {$2$};
        \node[state] at (6,0) (d) {$...$};
        \node[state] at (8,0) (n-2) {\tiny{$\R-2$}};
        \node[state] at (10,0) (n-1) {\tiny{$\R-1$}};
        \node[state] at (12,0) (n) {$\R$};
        
        \draw[every loop]
            (1) edge[bend right, auto=right] node {$1-p$} (0)
            (0) edge[loop above] node {$1$} (0)
            (n) edge[loop above] node {$1$} (n)
            (2) edge[bend right, auto=right] node {$1-p$} (1)
            (1) edge[bend right, auto=right] node {$p$} (2)
            (2) edge[bend right, auto=right] node {$p$} (d)
            (d) edge[bend right, auto=right] node {$p$} (n-2)
            (n-2) edge[bend right, auto=right] node {$p$} (n-1)
            (n-1) edge[bend right, auto=right] node {$p$} (n)
            (n-1) edge[bend right, auto=right] node {$p$} (n)
            (n) edge[bend right, auto=right] node {$1-p$} (n-1)
            (n-1) edge[bend right, auto=right] node {$1-p$} (n-2)
            (n-2) edge[bend right, auto=right] node {$1-p$} (d)
            (d) edge[bend right, auto=right] node {$1-p$} (2);
    \end{tikzpicture}
\captionof{figure}{Markov chain on infinite regular trees.} \label{fig:Markov_trees}
\end{figure}
    
\end{center}

Those familiar with the gambler's ruin problem will immediately realize that this Markov chain is the same as the gambler's ruin chain with a probability of the gambler winning a round given by  $p= t \left ( \frac{\Delta-1}{\Delta} \right)+r$ and the gambler losing a round given by $1-p = c+ \frac{t}{\Delta}$.

  It is well-known that the expected game length of the gambler's ruin problem is given by the equation $E(d) = d(\R-d)$ when $p=\frac{1}{2}$ \cite[p. 79]{D19}
and  $$E(d) = - \frac{\R}{1- 2p}\frac{\left( \left( \frac{1 - p}{p}\right)^d - 1\right)}{\left( \left( \frac{1 - p}{p}\right)^\R - 1\right)} + \frac{d}{1- 2p}$$ 
when $p\neq \frac{1}{2}$ \cite[p. 79]{D19}.

After finding a common denominator and factoring out $\frac{1}{1- 2p}$ we can rewrite $E(d)$ using our notation as

\begin{equation} E(d) = \left (\frac{1}{1 - 2p} \right)\left (d-\R \cdot p ^{\R-d} \cdot \left ( \frac{p^d - \left(1-p \right)^d}{p^\R - \left( 1-p \right)^\R} \right) \right) \label{excelp}  .\end{equation}
When $p< \frac{1}{2}$ (favoring cop success) then the right hand side of Equation \eqref{excelp} has a limit of $d$ as $\R \rightarrow \infty$. So if the cop never gives up, we get the formula
\begin{align*}
E(d) &= \frac{d}{1-2p} = \frac{d \Delta}{\Delta - 2\left(1-c\right)\left(\Delta-1\right)-2r}\
\\
 &= \frac{d \Delta}{\left(2c-1\right)\Delta + 2\left(1-c\right) - 2r}.
\end{align*}

 If  $c > \frac{1}{2}$, then as $\Delta \rightarrow \infty $ $$E(d) = \frac{d}{2c-1}.$$

If $R(d)$ and $C(d)$ represent the probability of the robber or cop winning respectively when starting the game from distance $d$ away from each other, then of course $R(d)+C(d)=1$. Formulas for $R(d)$ and $C(d)$ can easily be derived from well-known formulas modeling the gambler's ruin problem.  For instance,  when $p=\frac{1}{2}$, the probability of the robber escaping $R(d)$ is the same as the probability of the gambler's success in the fair gambler's ruin problem
$R(d) = \frac{d}{ \R} $ and when $p \not = \frac{1}{2}$,  we have 
$R(d) = \frac{1-\left(\frac{1-p}{p}\right)^d}{1-\left(\frac{1-p}{p}\right)^\R} $ is the same as the probability of the gambler's success in the unfair gambler's ruin problem \cite[Equation~1.2]{S09}.

Substituting $p = t \frac{\Delta-1}{\Delta} + r$ (transition probability to increase distance) and $\ds 1-p = c+\frac{t}{\Delta}$ (transition probability to decrease distance) we can find $R(d)$ to be
$$R(d)= \frac{1-\left(\frac{c+\frac{t}{\Delta}}{t \frac{\Delta-1}{\Delta} + r}\right)^d}{1-\left(\frac{c+\frac{t}{\Delta}}{t \frac{\Delta-1}{\Delta} +r}\right)^\R}.$$

\subsection{Numerical Examples}
As an example, we consider the game on an infinite regular tree of degree  $\Delta = 4$,  where the proportion of moves follows the distribution of $c=0.3$, $r=0.4$, $t=0.3$, and the cop calls off the chase if the robber reaches a distance of $\R=10$.  The following table gives the expected game time $E(d)$, and the probability of the robber or cop winning from each possible starting distance $d$.
\begin{center}
\begingroup
\setlength{\tabcolsep}{14pt}
    \begin{tabular}{|c|c|c|c|}\hline
        Distance & \multicolumn{3}{c|}{Measure}\\\hline
        $d$  & $E(d)$ & $R(d)$ & $C(d)$\\\hline
        $1$  & $12.10$ & $0.4024$ & $0.5976$\\\hline
        $2$  & $17.76$ & $0.6439$ & $0.3561$\\\hline
        $3$  & $19.55$ & $0.7888$ & $0.2112$\\\hline
        $4$  & $19.03$ & $0.8757$ & $0.1243$\\\hline
        $5$  & $17.11$ & $0.9278$ & $0.0722$\\\hline
        $6$  & $14.37$ & $0.9591$ & $0.0409$\\\hline
        $7$  & $11.12$ & $0.9779$ & $0.0221$\\\hline
        $8$  & $7.567$ & $0.9892$ & $0.0108$\\\hline
        $9$  & $3.838$ & $0.9959$ & $0.0041$\\\hline
    \end{tabular}
    \endgroup
\end{center}

We observe that if this game begins at $d=9$, then it has the lowest expected game time of any possible starting distance. This makes sense as $40\%$ of the time the robber will escape in the first move, or with probability $22.5\%$ either player will move drunkenly to allow the robber to escape in the first move; hence, the robber has a $66.5\%$ chance of escaping in the first move if the game begins at $d=9$. 

\section{General Matrix Method for Analyzing Markov Processes}  \label{sec:General_Method}

In this section we describe how to use various matrix equations to find the critical data points of the game.  These include, the probability of transitioning from one vertex to another in exactly $\M$ rounds, the probability the game lasts at least $\M$ rounds, and the expected game time. To calculate these values for a finite Markov chain we will use its probability matrix $P$, where $P_{i,j}$ is the probability of transitioning from state $i$ to state $j$ in one round of the Markov process.  We will use these matrix methods repeatedly to analyze various families of graphs in the subsequent sections of this paper.

\subsection{Calculating probability of moving from state $i$ to state $j$ in $\M$ rounds.}
The probability of starting in state $i$ and ending in state $j$ in exactly $\M$ moves is 
\begin{equation} P(i,j,\M) =  e_i \cdot P^{\M} \cdot e_j  \label{eq:transition} \end{equation}
where $e_i$ and $e_j$ are the row and column basis vectors respectively associated with state $i$ and state $j$.

For example, given the Markov chain in Figure \ref{fig:Markov_sample},
\begin{center}
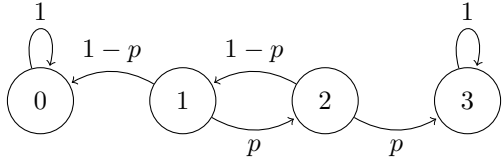
\begin{figure}[h] 
\begin{tikzpicture}
        \node[state]             (0) {0};
        \node[state, right=of 0] (1) {1};
        \node[state, right=of 1] (2) {2};
        \node[state, right=of 2] (3) {3};

        \draw[every loop]
            (1) edge[bend right, auto=right] node {$1-p$} (0)
            (0) edge[loop above] node {$1$} (0)
            (2) edge[bend right, auto=right] node {$1-p$} (1)
            (1) edge[bend right, auto=right] node {$p$} (2)
            (2) edge[bend right, auto=right] node {$p$} (3)
            (3) edge[loop above] node {$1$} (3);
     \end{tikzpicture}
\captionof{figure}{Sample for a Markov chain with four states.} \label{fig:Markov_sample}
\end{figure}
\end{center}
the probability matrix for this Markov chain is
\[P = 
\begin{pmatrix}
 1 & 0 & 0 & 0\\
 1-p & 0 & p & 0 \\
 0 & 1-p & 0 & p  \\
 0 & 0 & 0 & 1\\
\end{pmatrix}.
\]

 Some probabilities of going from one state to another are displayed in the table below
\begin{center}
\begin{tabular}{|c|c|c|c|c|}
\hline
    From state & to state & in $\M$ rounds & $P(i,j,\M)$ \\\hline
    $i = 3$ & $j = 0$ & $\M = 1$ & $P(3,0,1) = 0$ \\\hline
    $i = 3$ & $j = 0$ & $\M = 2$ & $P(3,0,2) = 0$ \\\hline
    $i = 3$ & $j = 0$ & $\M = 3$ & $P(3,0,3) = 0$ \\\hline
    $i = 2$ & $j = 0$ & $\M = 1$ & $P(2,0,1) = 0$ \\\hline
    $i = 2$ & $j = 0$ & $\M = 2$ & $\left(1-p\right)^2$ \\\hline
    $i = 2$ & $j = 0$ & $\M = 3$ & $\left(1-p\right)^2 \cdot 1 + \left(1-p\right)^3\cdot p$ \\\hline
    $i= 2$ & $j = 0$ & $\M = 4$ & $\left(1-p\right)^2 \cdot 1^2 + \left(1-p\right)^3\cdot p$ \\\hline
\end{tabular}

\end{center}
Note that in this example, state $0$ and state $3$ are both absorbing states, and that it is possible when calculating $P(i,a,\M)$ for an absorbing state $a$, that the robber reaches the absorbing state $a$ in $m<\M$ moves, and then sits there for the remaining $\M-m$ moves. 

\hypertarget{survival}{}
\subsection{Calculating the probability $\G_{\M}(d)$ that the game will last at least $\M$ rounds} \label{subsection:survival} 
To calculate the probability that the cop will still be actively chasing the robber through $\M$ rounds of this game with a starting distance $d$ between the players, we restrict attention to the matrix $T = P_{transient}$ modeling only the non-absorbing states of the Markov chain.  Note that in this case, $T$ is the matrix obtained from $P$ by removing the columns and rows associated with any absorbing states in $P$.
The probability of the game lasting at least $\M$ rounds of the game if the players start at distance $d$ from each other is then given by the following product of matrices \begin{equation} \G_{\M}(d) =  e_d \cdot T^{\M} \cdot \mathbbm{1}  \label{eq:survival} \end{equation} where $e_d$ denotes a standard basis row vector with $1$ in column $d$ and zero elsewhere, and $\ds \mathbbm{1}$ is the column vector with $1$ in each entry.

\hypertarget{expectation1}{}
\subsection{Calculating Expectation $E(d)$} \label{subsection:expectation}
As the number of rounds goes to infinity the likelihood that the game is still going on goes to zero. If we sum all the $\G_{\M}(d)$ probabilities for all $\M$ we can find the expected number of rounds the game should last. That is, \begin{align} \label{eq:expectation}  
E(d) & =  e_d \cdot \left [(I +T+T^2+T^3 + \cdots + T^{\M} +\cdots )\right ] \cdot \mathbbm{1}    \\ & =  e_d \cdot  \frac{1}{1-T} \cdot \mathbbm{1}  \\ 
& = e_d \cdot \left(I-T\right)^{-1} \cdot \mathbbm{1}
\end{align} 
where $I$ is the identity matrix.

\section{Cycle Graphs}\label{sec:Cycle Graphs}
The study of cycle graphs breaks down into two families of cases based on if the number of vertices $\n$ is even or odd. We start by analyzing the case when $\n$ is even, and then explain how the odd case differs slightly.  

For a cycle graph $C_{\n}$ with $\n$ nodes, where $\n$ is even we have the states where the robber and the cop are $0$ (cop catches the robber or the robber is tipsy and stumbles upon the cop herself), $1$, $2$, $3$, $\dots$ or $\frac{\n}{2}$ moves away from each other. 

We assume that the cop always moves until they are at distance 0 away. Also, if they are at distance $\frac{\n}{2}$ away, and the robber is sober, he will stay at the same location since it is the furthest from the cop. The Markov chain for a cycle graph with $\n$ nodes has  $\frac{\n}{2}$ states as shown in Figure \ref{fig:Markov_even_cycle}.
\begin{center}
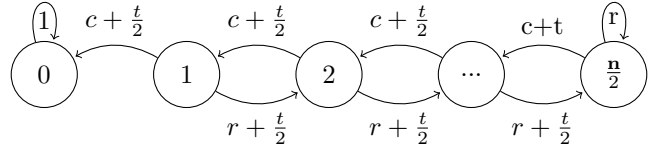
\begin{figure}[h] 
 \begin{tikzpicture}
        \node[state]             (0) {0};
        \node[state, right=of 0] (1) {1};
        \node[state, right=of 1] (2) {2};
        \node[state, right=of 2] (3) {...};
        \node[state, right=of 3] (4) {$\frac{\n}{2}$};

        \draw[every loop]
            (1) edge[bend right, auto=right] node {$c + \frac{t}{2}$} (0)
            (0) edge[loop above, auto=right] node {1} (1)
            
            (2) edge[bend right, auto=right] node {$c + \frac{t}{2}$} (1)
            (1) edge[bend right, auto=right] node {$r + \frac{t}{2}$} (2)
            (2) edge[bend right, auto=right] node {$r + \frac{t}{2}$} (3)
            (3) edge[bend right, auto=right] node {$c + \frac{t}{2}$} (2)
            (3) edge[bend right, auto=right] node {$r + \frac{t}{2}$} (4)
            (4) edge[loop above, auto=right] node {r} (4)
            (4) edge[bend right, auto=right ] node {c+t} (3);
    \end{tikzpicture}
\captionof{figure}{Markov chain cycle graph with $\n$ nodes, where $\n$ is even.} \label{fig:Markov_even_cycle}
\end{figure}
\end{center}
From this Markov chain we find the probability matrix $P$ to be a tridiagonal $ (\frac{\n}{2} + 1\times \frac{\n}{2} + 1)$ matrix of the following form: $P_{1,1} = 1$, $P_{\frac{\n}{2} + 1,\frac{\n}{2}} = c + t$, $P_{\frac{\n}{2} + 1,\frac{\n}{2} + 1} = r$ and the upper and lower diagonal entries are $P_{k+1,k} = c + \frac{t}{2}$ for $1 \leq k \leq \frac{\n}{2} - 1$, and  $P_{k+1,k+2} = r + \frac{t}{2}$ for $1 \leq k \leq \frac{\n}{2} - 1$, and finally $P_{i,j}=0$ otherwise.

  The transition matrix $T$ is derived by removing the first row and column corresponding to the absorbing state of $P$ . Hence, $T$ is an $(\frac{\n}{2} \times \frac{\n}{2})$ matrix with $T_{\frac{\n}{2},\frac{\n}{2} - 1} = c + t$, $T_{\frac{\n}{2},\frac{\n}{2}} = r$, $T_{k+1,k} = c + \frac{t}{2}$ for $1 \leq k \leq \frac{\n}{2} - 2$, $T_{k,k+1} = r + \frac{t}{2}$ for $1 \leq k \leq \frac{\n}{2} - 1$, and $T_{i,j}=0$ otherwise.

  When $\n$ is odd, the only difference in the Markov chain is that the probability of transitioning from state $\frac{\n}{2}$ to itself is $r+\frac{t}{2}$ and the probability of transitioning from $\frac{\n}{2}$ to $\frac{\n}{2} - 1$ is $c+\frac{t}{2}$.
Hence $P_{\frac{\n}{2} + 1,\frac{\n}{2}} = c + \frac{t}{2}$, $P_{\frac{\n}{2} + 1,\frac{\n}{2} + 1} = r + \frac{t}{2}$ , since the robber and cop move all the time. 

  We may now use our formulas involving $T$ to find the probability that the game lasts at least $\M$ rounds as in Subsection \ref{subsection:survival} and the expected game time Subsection \ref{subsection:expectation}.  It is also important to note that because our model does not guarantee turns alternate, and this graph is finite where the cop never calls off the chase, $C(d) = 1$ for all values of $d$, $c$, $r$, $t$, and $\n$.

\subsection{Numerical Examples}
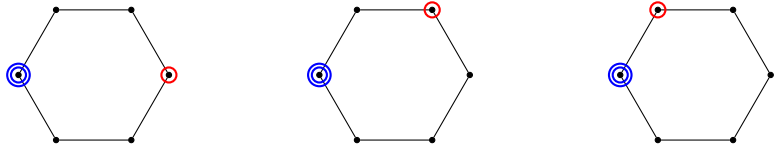
\begin{figure}[h]
\begin{tikzpicture}
\vertex[fill](a1) at (.5000000000, .8660254040) {};
\vertex[fill](a2) at (-.5000000000, .8660254040) {};
\vertex[fill](a3) at (-1., 0.) {};
\vertex[fill](a4) at (-.5000000000, -.8660254040) {};
\vertex[fill](a5) at (.5000000000, -.8660254040) {};
\vertex[fill](a6) at (1., 0.) {};
\draw (a1)--(a2) (a2)--(a3) (a3)--(a4) (a4)--(a5) (a5)--(a6) (a6)--(a1);
\draw[blue,thick] (-1,0) circle (.1cm);
\draw[blue,thick] (-1,0) circle (.15cm);
\draw[red,thick] (1,0) circle (.1cm);

\vertex[fill](a1) at (4.5000000000, .8660254040) {};
\vertex[fill](a2) at (3.5000000000, .8660254040) {};
\vertex[fill](a3) at (3., 0.) {};
\vertex[fill](a4) at (3.5000000000, -.8660254040) {};
\vertex[fill](a5) at (4.5000000000, -.8660254040) {};
\vertex[fill](a6) at (5., 0.) {};
\draw (a1)--(a2) (a2)--(a3) (a3)--(a4) (a4)--(a5) (a5)--(a6) (a6)--(a1);
\draw[blue,thick] (3,0) circle (.1cm);
\draw[blue,thick] (3,0) circle (.15cm);
\draw[red,thick] (4.5000000000, .8660254040) circle (.1cm);

\vertex[fill](a1) at (8.5000000000, .8660254040) {};
\vertex[fill](a2) at (7.5000000000, .8660254040) {};
\vertex[fill](a3) at (7., 0.) {};
\vertex[fill](a4) at (7.5000000000, -.8660254040) {};
\vertex[fill](a5) at (8.5000000000, -.8660254040) {};
\vertex[fill](a6) at (9., 0.) {};
\draw (a1)--(a2) (a2)--(a3) (a3)--(a4) (a4)--(a5) (a5)--(a6) (a6)--(a1);
\draw[blue,thick] (7,0) circle (.1cm);
\draw[blue,thick] (7,0) circle (.15cm);
\draw[red,thick] (7.5000000000, .8660254040) circle (.1cm);
\end{tikzpicture}
\captionof{figure}{Distance between cop (blue double circle) and robber (red single circle) on a cycle $C_6$ graph.}
 \label{fig:C_6}
\end{figure}

\hfill\\
\hypertarget{cycleFig}
For cycle graph with $\n=6$ nodes we have the states where the robber and the cop are 0 (cop catches the robber or the robber is tipsy and stumbles upon the cop himself), 1, 2 or 3 moves away from each other (Figure \ref{fig:C_6}). 

 We assume that the cop always moves until they are at distance 0 away. Also, if they are at distance 3 away and the robber is sober he will stay at the same location since it is the furthest from the cop. We can then create the Markov chain in Figure \ref{fig:Markov_C_6} for a cycle graph with 6 nodes.

\begin{center}
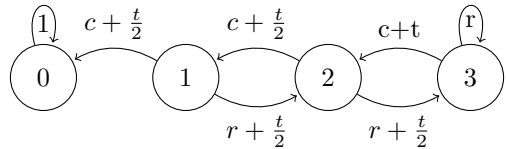
\begin{figure}[h] 
 \begin{tikzpicture}
        \node[state]             (0) {0};
        \node[state, right=of 0] (1) {1};
        \node[state, right=of 1] (2) {2};
        \node[state, right=of 2] (3) {3};

        \draw[every loop]
            (1) edge[bend right, auto=right] node {$c + \frac{t}{2}$} (0)
            (0) edge[loop above, auto=right] node {1} (1)
            (2) edge[bend right, auto=right] node {$c + \frac{t}{2}$} (1)
            (1) edge[bend right, auto=right] node {$r + \frac{t}{2}$} (2)
            (2) edge[bend right, auto=right] node {$r + \frac{t}{2}$} (3)
            (3) edge[loop above, auto=right] node {r} (3)
            (3) edge[bend right, auto=right ] node {c+t} (2);
    \end{tikzpicture}
\captionof{figure}{Markov chain on $C_6$ graph.} \label{fig:Markov_C_6}
\end{figure}
\end{center}
   
\begin{center}
 \begin{tabular}{| c | c | c | c | c |} 
 \hline
  & 0 & 1 & 2 & 3 \\ [0.5ex] 
 \hline
 0 & 1 & 0 & 0 & 0 \\ 
 \hline
 1 & $c + \frac{t}{2}$ & 0 & $r + \frac{t}{2}$ & 0 \\
 \hline
 2 & 0 & $c + \frac{t}{2}$ & 0 & $r + \frac{t}{2}$ \\
 \hline
 3 & 0 & 0 & $c + t$ & $r$ \\ [0.5ex] 
 \hline 
\end{tabular}
\end{center} 
We use a stochastic matrix, $P$ (below), to represent the transition probabilities of this system (rows and columns in this matrix are indexed by the possible states listed above, with the pre-transition state as the row and post-transition state as the column).
 Transition probability matrix
\[P = 
\begin{pmatrix}
 1 & 0 & 0 & 0\\
 c + \frac{t}{2} & 0 & r + \frac{t}{2} & 0 \\
 0 & c + \frac{t}{2} & 0 & r + \frac{t}{2}  \\
 0 & 0 & c + t & r\\
\end{pmatrix}
\]
This transition probability matrix $P$ can be restricted to a transient transition matrix $T$ with absorbing state 0 removed
\[T = 
\begin{pmatrix}
 0 & r + \frac{t}{2} & 0 \\
 c + \frac{t}{2} & 0 & r + \frac{t}{2}  \\
 0 & c + t & r\\
\end{pmatrix}
\]
Using our general matrix method as in Section \ref{sec:General_Method} we can solve the following numerical example.
\subsection{Numerical Examples}

 We assume the tipsy moves by either player account for half of all moves $t=0.5$ and $c+r=1-t=0.5$. The following tables gives the probability the cop will still be chasing the robber after $\M = 7$ rounds provided the cop and the robber start at distance $d = 1,2$ or $3$ as the percentage of the sober moves allocated to the cop and the robber varies.
\begin{center}
\begin{tabular}{|c|c|c|c|c|c|c|}\hline
    Measure & \multicolumn{6}{c|}{Proportion of Sober moves $c+r=0.5$} \\\hline 
    Robber  & $r=0$ & $r=0.1$ & $r=0.2$ & $r=0.3$ & $r=0.4$ & $r=0.5$ \\\hline
    Cop & $c=0.5$ & $c=0.4$ & $c=0.3$ & $c=0.2$ & $c=0.1$ & $c=0.0$ \\\hline
    $\G_{7}(1)$ & $0.020$ & $0.073$ & $0.169$ & $0.307$ & $0.473$ & $0.649$  \\\hline
    $\G_{7}(2)$ & $0.084$ & $0.195$ & $0.349$ & $0.526$ & $0.702$ & $0.849 $  \\\hline
    $\G_{7}(3)$ & $0.084$ & $0.222$ & $0.402$ & $0.59$ & $0.759$ & $0.887$  \\\hline
    $E(1)$ & $1.89$ & $2.69$ & $4.14$ & $7.07$ & $13.91$ & $34$  \\\hline
    $E(2)$ & $3.56$ & $4.83$ & $6.98$ & $11.04$ & $19.86$ & $44$  \\\hline
    $E(3)$ & $4.56$ & $5.94$ & $8.23$ & $12.47$ & $21.53$ & $46$  \\\hline
\end{tabular}
\end{center}
The table shows the game will last longest, and the chase has the largest probability of still continuing after $7$ rounds, if the starting distance is $3$ and the robber takes the maximum percentage of sober moves possible.
The accompanying CoCalc code for these computations is included in Appendix \ref{sec:code_Petersen}, so the interested reader can adapt these calculations to model the game for any values of $\M, r, c, t$ they choose.


\section{Petersen Graph}\label{sec:Petersen}

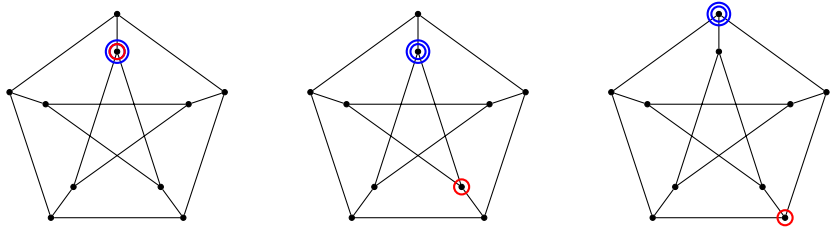
\begin{figure}[h]

\begin{tikzpicture}
\vertex[fill](a1) at (0,1) {};
\vertex[fill](o1) at (0,1.5) {};
\vertex[fill](a2) at (-.95,.30) {};
\vertex[fill](o2) at (-1.43,0.46) {};
\vertex[fill](a3) at (.95,.30) {};
\vertex[fill](o3) at (1.43,0.46) {};
\vertex[fill](a4) at (-.58,-.80) {};
\vertex[fill](o4) at (-0.88,-1.21) {};
\vertex[fill](a5) at (.58,-.80) {};
\vertex[fill](o5) at (0.88,-1.21) {};
\draw (a1)--(a4) (a4)--(a3) (a3)--(a2) (a5)--(a2) (a5)--(a1);
\draw (o1)--(o2) (o2)--(o4) (o4)--(o5) (o5)--(o3) (o3)--(o1);
\draw (a1)--(o1) (a2)--(o2) (a3)--(o3) (a4)--(o4) (a5)--(o5);
\draw[blue,thick] (a1) circle (.1cm);
\draw[blue,thick] (a1) circle (.15cm);
\draw[red,thick] (a1) circle (.1cm);

\vertex[fill](a1) at (4,1) {};
\vertex[fill](o1) at (4,1.5) {};
\vertex[fill](a2) at (3.05,.30) {};
\vertex[fill](o2) at (2.57,0.46) {};
\vertex[fill](a3) at (4.95,.30) {};
\vertex[fill](o3) at (5.43,0.46) {};
\vertex[fill](a4) at (3.42,-.80) {};
\vertex[fill](o4) at (3.12,-1.21) {};
\vertex[fill](a5) at (4.58,-.80) {};
\vertex[fill](o5) at (4.88,-1.21) {};
\draw (a1)--(a4) (a4)--(a3) (a3)--(a2) (a5)--(a2) (a5)--(a1);
\draw (o1)--(o2) (o2)--(o4) (o4)--(o5) (o5)--(o3) (o3)--(o1);
\draw (a1)--(o1) (a2)--(o2) (a3)--(o3) (a4)--(o4) (a5)--(o5);
\draw[blue,thick] (a1) circle (.1cm);
\draw[blue,thick] (a1) circle (.15cm);
\draw[red,thick] (a5) circle (.1cm);

\vertex[fill](a1) at (8,1) {};
\vertex[fill](o1) at (8,1.5) {};
\vertex[fill](a2) at (7.05,.30) {};
\vertex[fill](o2) at (6.57,0.46) {};
\vertex[fill](a3) at (8.95,.30) {};
\vertex[fill](o3) at (9.43,0.46) {};
\vertex[fill](a4) at (7.42,-.80) {};
\vertex[fill](o4) at (7.12,-1.21) {};
\vertex[fill](a5) at (8.58,-.80) {};
\vertex[fill](o5) at (8.88,-1.21) {};
\draw (a1)--(a4) (a4)--(a3) (a3)--(a2) (a5)--(a2) (a5)--(a1);
\draw (o1)--(o2) (o2)--(o4) (o4)--(o5) (o5)--(o3) (o3)--(o1);
\draw (a1)--(o1) (a2)--(o2) (a3)--(o3) (a4)--(o4) (a5)--(o5);
\draw[blue,thick] (o1) circle (.1cm);
\draw[blue,thick] (o1) circle (.15cm);
\draw[red,thick] (o5) circle (.1cm);
\end{tikzpicture}
\captionof{figure}{Distance between cop (blue double circle) and robber (red single circle) on a Petersen graph.}
\label{fig:P}
\end{figure}

 The Petersen graph is a vertex-transitive graph where the cop and robber can only be  $0$, $1$, or $2$ moves away from each other as (Figure \ref{fig:P}).  

 We assume that the cop moves until eventually she captures the robber. Also, if robber is distance 2 away from the cop, and the robber is sober, he will stay at the same location since it is the furthest from the cop. Hence the Markov chain for the Petersen graph is as depicted in Figure \ref{fig:Markov_P}.
\begin{center}
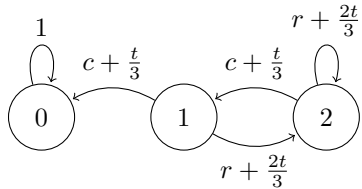
\begin{figure}[h] 
 \begin{tikzpicture}
        \node[state]             (0) {0};
        \node[state, right=of 0] (1) {1};
        \node[state, right=of 1] (2) {2};
        
        \draw[every loop]
            (1) edge[bend right, auto=right] node {$c + \frac{t}{3}$} (0)
            (0) edge[loop above] node {1} (1)
            
            (2) edge[bend right, auto=right] node {$c + \frac{t}{3}$} (1)
            (1) edge[bend right, auto=right] node {$r + \frac{2t}{3}$} (2)
            (2) edge[loop above] node {$r + \frac{2t}{3}$} (2);
    \end{tikzpicture}
\captionof{figure}{Markov chain for Petersen graph.} \label{fig:Markov_P}
\end{figure}
\end{center}   
    
\begin{center}
 \begin{tabular}{| c | c | c | c |} 
 \hline
  & 0 & 1 & 2 \\ [0.5ex] 
 \hline
 0 & 1 & 0 & 0 \\ 
 \hline
 1 & $c + \frac{t}{3}$ & 0 & $r + \frac{2t}{3}$\\
 \hline
 2 & 0 & $c + \frac{t}{3}$ & $r + \frac{2t}{3}$ \\ [0.5ex] 
 \hline 
\end{tabular}
\end{center} 
We use a stochastic matrix $P$ to represent the transition probabilities of this system. The rows and columns in this matrix are indexed by the possible states listed above, with the pre-transition state as the row and the post-transition state as the column.
\[P = 
\begin{pmatrix}
 1 & 0 & 0 \\
 c + \frac{t}{3} & 0 & r + \frac{2t}{3} \\
 0 & c + \frac{t}{3} & r + \frac{2t}{3} \\
\end{pmatrix}
\]
The states where there is a distance between the cop and the robber do not contribute to the survival average so state 0 can be ignored. The initial state and transition matrix can be reduced to a transition matrix with absorbing state 0 removed
\[T = 
\begin{pmatrix}
 0 & r + \frac{2t}{3} \\
 c + \frac{t}{3} & r + \frac{2t}{3}  \\
\end{pmatrix}
\]

 We may now use our formulas involving $T$ to find the probability that the game will last at least $\M$ rounds as in Subsection \ref{subsection:survival} and the expected number of rounds as in Subsection \ref{subsection:expectation}.

\subsection{Numerical Examples}\hfill\\
 In the table below, we assume the tipsy moves by either player account for half of all moves $t=0.5$ and $c+r=1-t=0.5$.  We then calculate the probability the cop will still be chasing the robber after $\M = 7$ rounds based on their starting distances $d = 1$ or $2$ and on the percentage of the sober moves allocated to the cop and the robber.
\begin{center}
\begin{tabular}{|c|c|c|c|c|c|c|}\hline
    Measure & \multicolumn{6}{c|}{Proportion of Sober moves $c+r=0.5$} \\\hline 
    Robber  & $r=0$ & $r=0.1$ & $r=0.2$ & $r=0.3$ & $r=0.4$ & $r=0.5$ \\\hline
    Cop & $c=0.5$ & $c=0.4$ & $c=0.3$ & $c=0.2$ & $c=0.1$ & $c=0.0$ \\\hline
    $\G_{7}(1)$ & $0.039$ & $0.1$ & $0.204$ & $0.352$ & $0.536$ & $0.735$  \\\hline
    $\G_{7}(2)$ & $0.078$ & $0.175$ & $0.318$ & $0.496$ & $0.687$ & $0.86$  \\\hline
    $E(1)$ & $2.25$ & $3.11$ & $4.59$ & $7.44$ & $14.06$ & $40$  \\\hline
    $E(2)$ & $3.75$ & $4.87$ & $6.73$ & $10.17$ & $17.81$ & $42$  \\\hline
\end{tabular}
\end{center}

We have published the accompanying CoCalc code for these computations in Appendix \ref{sec:code_Petersen}, so the interested reader can change these calculations to model the game for any values of $\M, r, c, t$ that they choose.

\section{Friendship Graphs}\label{sec:Friendship Graphs}
As friendship graphs are not vertex-transitive, we cannot simply model the game on them with a Markov chain where each state is determined by the distance between the cop and robber, and we must treat tipsy cop moves and tipsy robber moves separately as the transition probabilities from state to state depend on who is moving.  Hence, in this section the spinner we use has four distinct possible outcomes satisfying $r+c + \tr+\tc=1$. 

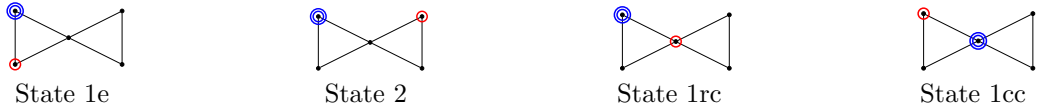
\begin{figure}[h]
\begin{minipage}{\textwidth}
\begin{minipage}[b]{0.2\textwidth}
\centering
\resizebox{.5\textwidth}{!}{
\begin{tikzpicture}
\vertex[fill](a0) at (0,0) {};
\vertex[fill](b1) at (-1,.5) {};
\vertex[fill](b2) at (-1,-.5) {};
\vertex[fill](c1) at (1,.5) {};
\vertex[fill](c2) at (1,-.5) {};
\draw (b1)--(a0) (b1)--(b2) (b2)--(a0);
\draw (c1)--(a0) (c1)--(c2) (c2)--(a0);
\draw[blue,thick] (-1,.5) circle (.1cm);
\draw[blue,thick] (-1,.5) circle (.15cm);
\draw[red,thick] (-1,-.5) circle (.1cm);
\end{tikzpicture}}\\
State 1e
\end{minipage}
\hspace{5mm}
\begin{minipage}[b]{0.2\textwidth}
\centering
\resizebox{.5\textwidth}{!}{
\begin{tikzpicture}
\vertex[fill](a0) at (0,0) {};
\vertex[fill](b1) at (-1,.5) {};
\vertex[fill](b2) at (-1,-.5) {};
\vertex[fill](c1) at (1,.5) {};
\vertex[fill](c2) at (1,-.5) {};
\draw (b1)--(a0) (b1)--(b2) (b2)--(a0);
\draw (c1)--(a0) (c1)--(c2) (c2)--(a0);
\draw[blue,thick] (-1,.5) circle (.1cm);
\draw[blue,thick] (-1,.5) circle (.15cm);
\draw[red,thick] (1,.5) circle (.1cm);
\end{tikzpicture}}\\
State 2
\end{minipage}
\hspace{5mm}
\begin{minipage}[b]{0.2\textwidth}
\centering
\resizebox{.5\textwidth}{!}{
\begin{tikzpicture}
\vertex[fill](a0) at (0,0) {};
\vertex[fill](b1) at (-1,.5) {};
\vertex[fill](b2) at (-1,-.5) {};
\vertex[fill](c1) at (1,.5) {};
\vertex[fill](c2) at (1,-.5) {};
\draw (b1)--(a0) (b1)--(b2) (b2)--(a0);
\draw (c1)--(a0) (c1)--(c2) (c2)--(a0);
\draw[blue,thick] (-1,.5) circle (.1cm);
\draw[blue,thick] (-1,.5) circle (.15cm);
\draw[red,thick] (0,0) circle (.1cm);
\end{tikzpicture}}\\
State 1rc
\end{minipage}
\hspace{5mm}
\begin{minipage}[b]{0.2\textwidth}
\centering
\resizebox{.5\textwidth}{!}{
\begin{tikzpicture}
\vertex[fill](a0) at (0,0) {};
\vertex[fill](b1) at (-1,.5) {};
\vertex[fill](b2) at (-1,-.5) {};
\vertex[fill](c1) at (1,.5) {};
\vertex[fill](c2) at (1,-.5) {};
\draw (b1)--(a0) (b1)--(b2) (b2)--(a0);
\draw (c1)--(a0) (c1)--(c2) (c2)--(a0);
\draw[red,thick] (-1,.5) circle (.1cm);
\draw[blue,thick] (0,0) circle (.1cm);
\draw[blue,thick] (0,0) circle (.15cm);
\end{tikzpicture}}\\
State 1cc
\end{minipage}
\end{minipage}
\captionof{figure}{All possible states of Friendship Graphs of $n=2$ triangles.}
\label{fig:Friendship3}

\end{figure}

The following notation will be used to model the tipsy cop and tipsy robber game on friendship graphs of $n$ triangles. The states $1e$, $1cc$, and $1rc$ all refer to the cop and robber being distance $1$ away from each other but, $1e$ means both players are on the same outer edge, $1cc$ means the cop is in the center, and $1rc$  means the robber is in the center as depicted in Figure \ref{fig:Friendship3}.

\begin{figure}[h]
\centering
\begin{tikzpicture}
\vertex[fill](a1) at (.75,.75) {};
\vertex[fill](a2) at (-.75,.75) {};
\vertex[fill](a3) at (-.75,-.75) {};
\vertex[fill](a4) at (.75,-.75) {};
\vertex[fill](a5) at (0,0) {};
\draw (a1)--(a2) ;
\draw (a2)--(a4);
\draw (a3)--(a1);
\draw (a3)--(a4) ;
\draw[blue,thick] (a1) circle (.1cm);
\draw[blue,thick] (a1) circle (.15cm);
\draw[red,thick] (a3) circle (.1cm);
\end{tikzpicture}
\hspace{1cm}
\begin{tikzpicture}
\vertex[fill](a1) at (.5000000000, .8660254040) {};
\vertex[fill](a2) at (-.5000000000, .8660254040) {};
\vertex[fill](a3) at (-1., 0.) {};
\vertex[fill](a4) at (-.5000000000, -.8660254040) {};
\vertex[fill](a5) at (.5000000000, -.8660254040) {};
\vertex[fill](a6) at (1., 0.) {};
\vertex[fill](a0) at (0., 0.) {};
\draw (a1)--(a2) ;
\draw (a0)--(a2) ;
\draw (a3)--(a4) ;
\draw (a0)--(a4) ;
\draw (a3)--(a0) ;
\draw (a0)--(a1) ;
\draw (a5)--(a6);
\draw (a5)--(a0);
\draw (a0)--(a6);
\draw[blue,thick] (-1,0) circle (.1cm);
\draw[blue,thick] (-1,0) circle (.15cm);
\draw[red,thick] (1,0) circle (.1cm);
\end{tikzpicture}
\hspace{1cm}
\begin{tikzpicture}
\vertex[fill](a0) at (0., 0.) {};
\vertex[fill](a1) at (.707107, .707107) {};
\vertex[fill](a2) at (0,1) {};
\vertex[fill](a3) at (-0.707107,0.707107) {};
\vertex[fill](a4) at (-1,0) {};
\vertex[fill](a5) at (-0.707107,-0.707107) {};
\vertex[fill](a6) at (0,-1) {};
\vertex[fill](a7) at (0.707107,-0.707107) {};
\vertex[fill](a8) at (1,0) {};
\draw (a0)--(a1);
\draw (a0)--(a2);
\draw (a0)--(a3);
\draw (a0)--(a4);
\draw (a0)--(a5);
\draw (a0)--(a6);
\draw (a0)--(a7);
\draw (a0)--(a8);
\draw (a1)--(a2);
\draw (a3)--(a4);
\draw (a5)--(a6);
\draw (a7)--(a8);
\draw[blue,thick] (-0.707107,-0.707107) circle (.1cm);
\draw[blue,thick] (-0.707107,-0.707107) circle (.15cm);
\draw[red,thick] (0.707107,0.707107) circle (.1cm);
\end{tikzpicture}
\captionof{figure}{Examples of friendship graphs with $n=2,3,4$ triangles with cop (blue) and robber (red).}
    \label{fig:Friendshipn}
\end{figure}
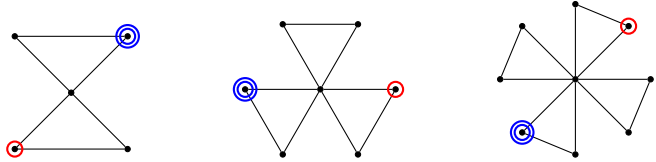

  Friendship graphs with $n=2,3,$ and $4$ triangles are shown in Figure \ref{fig:Friendshipn}. 
  Since the cop never calls off the chase on this finite graph, the cop will win eventually $C(d) = 1$ even if the robber is completely sober throughout.

The Markov chain modeling the game on a friendship graph with $n$ triangles is shown in Figure \ref{fig:markovchainfriendshipgraphs}.

\begin{center}
\begin{figure}[h] 
\begin{tikzpicture}
        \node[state] at (0,0)       (0) {0};
        \node[state] at (3,0)     (1e) {1e};
        \node[state] at (5,2)   (1cc) {1cc};
        \node[state] at (5,-2)  (1rc) {1rc};
        \node[state] at (8,0)       (2) {2};
        
        \draw[every loop]
            (0) edge[loop left] node {$1$} (0)
            (1e) edge[ auto=right] node {$c+\frac{\tc}{2}+\frac{\tr}{2}$} (0)
            (1cc) edge[bend right, auto=right] node {$c+\frac{\tc}{2n}+\frac{\tr}{2}$} (0)
            (1rc) edge[bend left, auto=left] node {$c+\frac{\tr}{2n}+ \frac{\tc}{2}$} (0)
            (2) edge[bend left, auto=right] node {$c+\frac{\tc}{2}$} (1cc)
            (1e) edge[bend right, auto=left] node {$\frac{\tc}{2}$} (1cc)
            (1cc) edge[loop above] node {$r+\frac{\tr}{2}$} (1cc)
            (2) edge[loop right] node {$r+\frac{\tc}{2}+\frac{\tr}{2}$} (2)
            (1cc) edge[bend left, auto=left] node {$\tc \cdot\frac{n-1}{n}$} (2)
            (2) edge[bend right, auto=left] node {$\frac{\tr}{2}$} (1rc)
            (1rc) edge[bend right, auto=right] node {$r+\tr\cdot \frac{n-1}{n}$} (2)
            (1rc) edge[loop below] node {$\frac{\tc}{2}$} (1rc)
            (1rc) edge[bend left, auto=left] node {$\frac{\tr}{2n}$} (1e)
            (1e) edge[bend left, auto=left] node {$r+\frac{\tr}{2}$} (1rc)
            (1cc) edge[bend right, auto=right] node {$\frac{\tc}{2n}$} (1e);
     \end{tikzpicture}
     \captionof{figure}{Markov chain for tipsy cop and tipsy robber game on friendship graphs of $n$ cliques}
    \label{fig:markovchainfriendshipgraphs}
\end{figure}
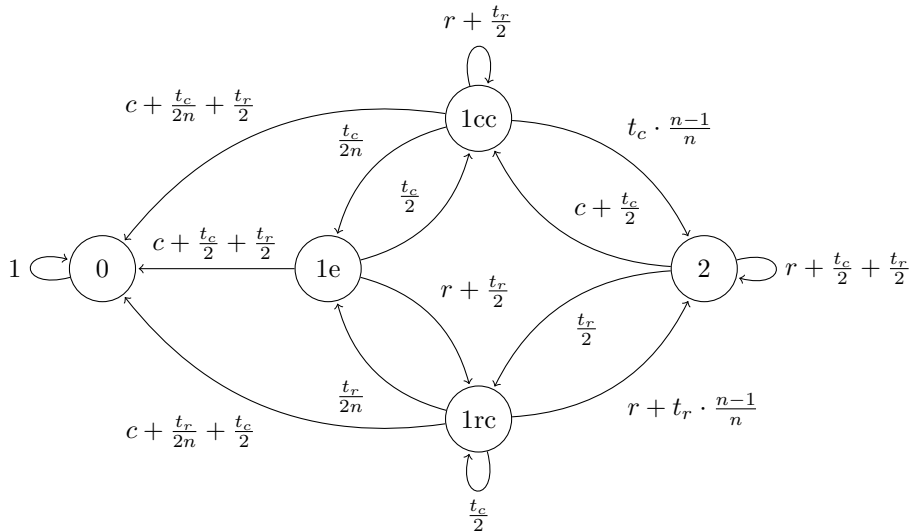
\end{center}

\begin{center}
    \begin{tabular}{|c|c|c|c|c|c|}
    \hline
    From \textbackslash To  & $2$ & $1cc$ & $1rc$ & $1e$ & $0$ \\\hline
         $2$  & $r+\frac{\tc}{2}+\frac{\tr}{2}$ & $c+\frac{\tc}{2}$ & $\frac{\tr}{2}$ & $0$ & $0$\\\hline
        $1cc$ &$\tc \cdot \frac{n-1}{n}$ & $r+\frac{\tr}{2}$ & $0$ & $\frac{\tc}{2n}$ & $c+\frac{\tc}{2n}+\frac{\tr}{2}$\\\hline
        $1rc$ & $r+\tr\cdot \frac{n-1}{n}$ & $0$ & $\frac{\tc}{2}$ & $\frac{\tr}{2n}$ & $c+\frac{\tr}{2n}+\frac{\tc}{2}$\\\hline
        $1e$  & $0$ & $\frac{\tc}{2}$ & $r+\frac{\tr}{2}$ & $0$ & $c+ \frac{\tc}{2}+\frac{\tr}{2}$\\\hline
        $0$   & $0$ & $0$ & $0$ & $0$ & $1$\\\hline
    \end{tabular}
\end{center}

From this Markov chain we derive the transition probability matrix  $P$ and transition matrix with absorbing state 0 removed $T$
\[P = 
\begin{pmatrix}
         r+\frac{\tc}{2}+\frac{\tr}{2} & c+\frac{\tc}{2} & \frac{\tr}{2} & 0 & 0\\
        \tc \cdot \frac{n-1}{n} & r+\frac{\tr}{2} & 0 & \frac{\tc}{2n} & c+\frac{\tc}{2n}+\frac{\tr}{2}\\
         r+\tr\cdot \frac{n-1}{n} & 0 & \frac{\tc}{2} & \frac{\tr}{2n} & c+\frac{\tr}{2n}+\frac{\tc}{2}\\
         0 & \frac{\tc}{2} & r+\frac{\tr}{2} & 0 & c+ \frac{\tc}{2}+\frac{\tr}{2}\\
         0 & 0 & 0 & 0 & 1 \\
\end{pmatrix} \quad
T= \begin{pmatrix}
         r+\frac{\tc}{2}+\frac{\tr}{2} & c+\frac{\tc}{2} & \frac{\tr}{2} & 0 \\
        \tc \cdot \frac{n-1}{n} & r+\frac{\tr}{2} & 0 & \frac{\tc}{2n} \\
         r+\tr\cdot \frac{n-1}{n} & 0 & \frac{\tc}{2} & \frac{\tr}{2n} \\
         0 & \frac{\tc}{2} & r+\frac{\tr}{2} & 0 & \\
\end{pmatrix}
\]

\subsection{Sample calculations} 
In this subsection we compute the expected game times and probability of the game lasting through $\M=10$ rounds on a friendship graph with $n=5$ cliques.  In these computations we assume that the cop and robber each get $50\%$ of the moves so $r+\tr=c+\tc=0.5$, but we vary the tipsiness of the cop and robber.  The code for these computations is available in Appendix \ref{sec:code_friendship}.

\begin{center}
\begin{tabular}{|c|c|c|c|c|c|c|c|c|c|}\hline
    Measure & \multicolumn{9}{c|}{Proportion of Tipsy moves} \\\hline 
    Robber  & \multicolumn{3}{c|}{$\tr=0.1$} & \multicolumn{3}{c|}{$\tr=0.25$} & \multicolumn{3}{c|}{$\tr=0.4$}\\\hline
    Cop& $\tc=0.1$ & $\tc=0.25$ & $\tc=0.4$ & $\tc=0.1$ & $\tc=0.25$ & $\tc=0.4$ & $\tc=0.1$ & $\tc=0.25$ & $\tc=0.4$\\\hline
    $\G_{10}(2)$ & $0.0448$ & $0.1836$ & $0.4300$ & $0.0252$ & $0.1129$ & $0.2836$ & $0.0143$ & $0.0691$ & $0.1862$ \\\hline
    $\G_{10}(1cc)$ & $0.0156$ & $0.1043$ & $0.3198$ & $0.0076$ & $0.0584$ & $0.1948$ & $0.0037$ & $0.0326$& $0.1183$ \\\hline
    $\G_{10}(1rc)$ & $0.0337$ & $0.1275$ & $0.2949$ & $0.0192$ & $0.0791$ & $0.1960$ & $0.0109$ & $0.0486$ &$0.1291$ \\\hline
    $\G_{10}(1e)$ & $0.0229$ & $0.0856$ & $0.2167$ & $0.0116$ & $0.0473$ & $0.1297$ & $0.0056$ & $0.0251$ & $0.0752$ \\\hline
    $E(2)$ & $4.628$ & $6.914$ & $12.09$ & $4.093$ & $5.697$ & $8.668$ & $3.689$ & $4.893$ & $6.874$ \\\hline
    $E(1cc)$ & $2.540$ & $4.505$ & $9.344$ & $2.160$ & $3.543$ & $6.333$ & $1.878$ & $2.914$ & $4.769$ \\\hline
    $E(1rc)$ & $3.419$ & $4.979$ & $8.590$ & $3.051$ & $4.159$ & $6.270$ & $2.765$ & $3.603$ & $5.032$ \\\hline
    $E(1e)$ & $2.665$ & $3.804$ & $6.735$ & $2.252$ & $3.002$ & $4.618$ & $1.923$ & $2.445$ & $3.463$ \\\hline
\end{tabular}
\end{center}

\section{Toroidal Grids} \label{sec:Toroidal}
A toroidal grid is the Cartesian product of two cycle graphs $C_m \square C_n$ and represents a grid that can be embedded on the surface of a torus, as shown in Figure \ref{fig:Toroidal Grid}.  Since our model does not guarantee turns alternate, and the cop never calls off the chase, the cop will eventually win $C(d) = 1$ on any toroidal grid.

\begin{center}
    \includegraphics[scale=0.5]{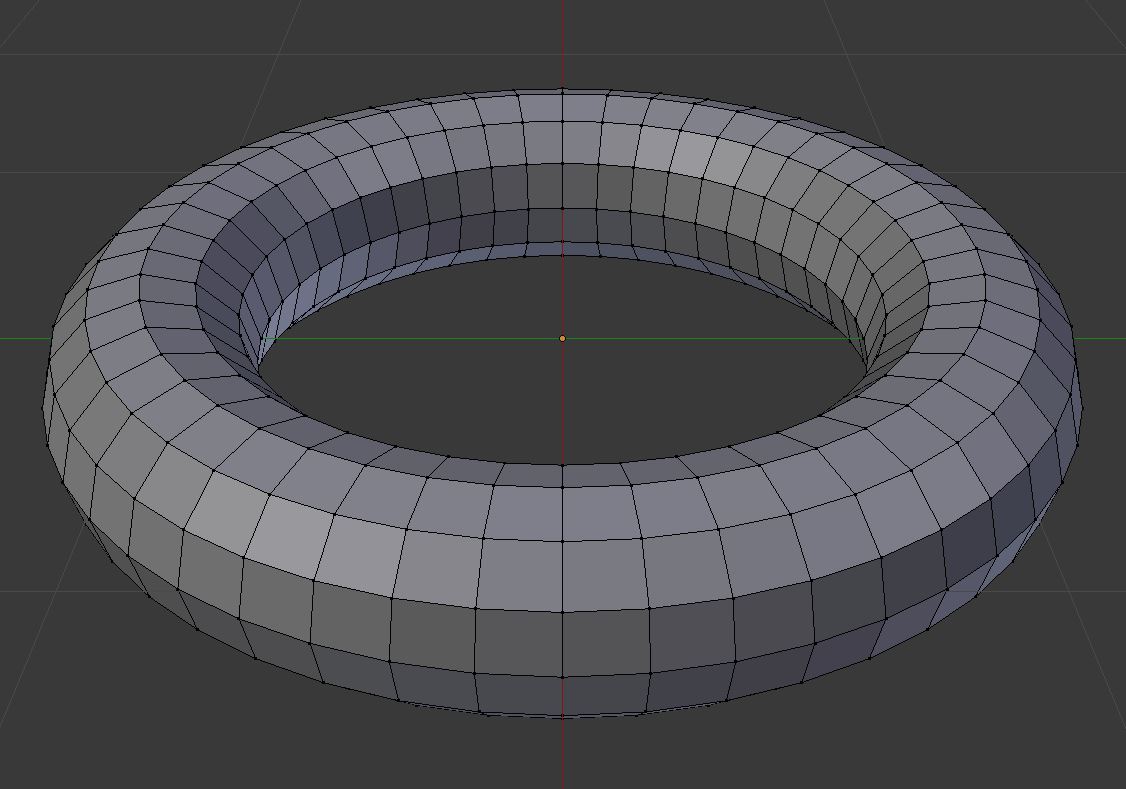}
    \captionof{figure}{Toroidal Grid. Generated using Blender open-source imaging software.}
    \label{fig:Toroidal Grid}
\end{center}

 The Markov chain for our model on a $7 \times 7$ toroidal grids is shown in Figure \ref{fig:Markov_toroidal7}. The number of states in our Markov chain model grows considerably as the number of nodes in $C_m \times C_n$ grows.
\begin{center}
\begin{figure}[h] 
\begin{tikzpicture}
        \node[state] at (0,0)       (0) {0};
        \node[state] at (2,0)     (10) {1,0};
        \node[state] at (4,0)  (20) {2,0};
        \node[state] at (6,0)  (30) {3,0};
        \node[state] at (2,2)  (11) {1,1};
        \node[state] at (4,2)  (21) {2,1};
        \node[state] at (6,2)  (31) {3,1};
        \node[state] at (4,4)  (22) {2,2};
        \node[state] at (6,4)  (32) {3,2};
        \node[state] at (6,6)  (33) {3,3};
        
    \draw[every loop]    
    
(0) edge[loop left] node {$1$} (0)
(10) edge[bend left, auto = left] node {$c + \frac{t}{4}$}(0)
(10) edge[bend right, auto = left] node {$\frac{t}{2}$} (11)
(10) edge[bend left, auto = right] node {$r+\frac{t}{4}$} (20)
(20) edge[bend left, auto = left] node {$c+\frac{t}{4}$} (10)
(30) edge[bend left, auto = left] node {$c+\frac{t}{4}$} (20)
(11) edge[bend left, auto = right] node {$r+\frac{t}{2}$} (21)
(20) edge[bend right, auto = left] node {$\frac{t}{2}$} (21)
(21) edge[bend right, auto = left] node {$\frac{t}{4}$} (22)
(22) edge[bend right, auto = right] node {$c+\frac{t}{2}$} (21)
(21) edge[bend left, auto = left] node {$c+\frac{t}{4}$} (11)
(31) edge[bend left, auto = left] node {$c+\frac{t}{4}$} (21)
(21) edge[bend right, auto = left] node {$\frac{t}{4}$} (20)
(11) edge[bend right, auto = right] node {$c+\frac{t}{2}$} (10)
(32) edge[bend left, auto = left] node {$c+\frac{t}{4}$} (22)
(30) edge[loop right] node {$\frac{t}{4}$} (30)
(31) edge[loop right] node {$\frac{t}{4}$} (31)
(32) edge[loop right] node {$\frac{t}{4}$} (32)
(33) edge[loop right] node {$r+\frac{t}{2}$} (33)
(20) edge[bend left, auto = right] node {$r+\frac{t}{4}$} (30)
(21) edge[bend left, auto = right] node {$r+\frac{t}{4}$} (31)
(22) edge[bend left, auto = right] node {$r+\frac{t}{2}$} (32)
(30) edge[bend right, auto = right] node {$r+\frac{t}{2}$} (31)
(31) edge[bend right, auto = right] node {$r+\frac{t}{4}$} (32)
(32) edge[bend right, auto = right] node {$r+\frac{t}{4}$} (33)
(33) edge[bend right, auto = right] node {$c+\frac{t}{2}$} (32)
(32) edge[bend right, auto = left] node {$\frac{t}{4}$} (31)
(31) edge[bend right, auto = left] node {$\frac{t}{4}$} (30);

     \end{tikzpicture}\\
   \captionof{figure}{Markov chain for a $7 \times 7$ toroidal grid.}  
\label{fig:Markov_toroidal7}
\end{figure}
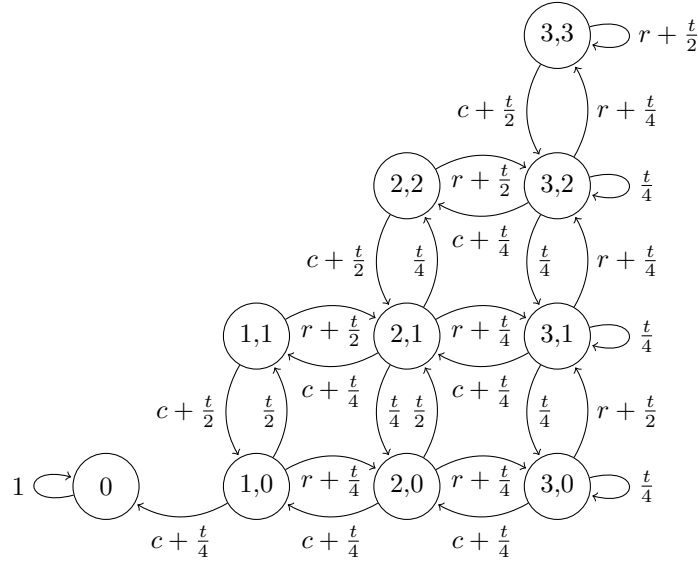
\end{center}

\begin{center}
    \begin{tabular}{|c|c|c|c|c|c|c|c|c|c|c|}
    \hline
        From \textbackslash To & (3,3)&(3,2)&(3,1)&(3,0)&(2,2)&(2,1)&(2,0)&(1,1)&(1,0)&0 \\\hline
        (3,3) & $r+\frac{t}{2}$&$c+\frac{t}{2}$&$0$&$0$&$0$&$0$&$0$&$0$&$0$&$0$\\\hline
        (3,2)&$r+\frac{t}{4}$&$\frac{t}{4}$&$\frac{t}{4}$&$0$&$c+\frac{t}{4}$&$0$&$0$&$0$&$0$&$0$\\\hline
        (3,1)&$0$&$r+\frac{t}{4}$&$\frac{t}{4}$&$\frac{t}{4}$&$0$&$c+\frac{t}{4}$&$0$&$0$&$0$&$0$\\\hline
        (3,0)&$0$&$0$&$r+\frac{t}{2}$&$\frac{t}{4}$&$0$&$0$&$c+\frac{t}{4}$&$0$&$0$&$0$\\\hline
        (2,2)&$0$&$r+\frac{t}{2}$&$0$&$0$&$0$&$c+\frac{t}{2}$&$0$&$0$&$0$&$0$\\\hline
        (2,1)&$0$&$0$&$r+\frac{t}{4}$&$0$&$\frac{t}{4}$&$0$&$\frac{t}{4}$&$c+\frac{t}{4}$&$0$&$0$\\\hline
        (2,0)&$0$&$0$&$0$&$r+\frac{t}{4}$&$0$&$\frac{t}{2}$&$0$&$0$&$c+\frac{t}{4}$&$0$\\\hline
        (1,1)&$0$&$0$&$0$&$0$&$0$&$r+\frac{t}{2}$&$0$&$0$&$c+\frac{t}{2}$&$0$\\\hline
        (1,0)&$0$&$0$&$0$&$0$&$0$&$0$&$r+\frac{t}{4}$& $\frac{t}{2}$&$0$&$c+\frac{t}{4}$\\\hline
        0&$0$&$0$&$0$&$0$&$0$&$0$&$0$&$0$&$0$&$1$\\\hline
    \end{tabular}
\end{center}
Therefore, \[ P = \begin{pmatrix}
        r+\frac{t}{2}&c+\frac{t}{2}&0&0&0&0&0&0&0&0\\
        r+\frac{t}{4}&\frac{t}{4}&\frac{t}{4}&0&c+\frac{t}{4}&0&0&0&0&0\\
        0&r+\frac{t}{4}&\frac{t}{4}&\frac{t}{4}&0&c+\frac{t}{4}&0&0&0&0\\
        0&0&r+\frac{t}{2}&\frac{t}{4}&0&0&c+\frac{t}{4}&0&0&0\\
        0&r+\frac{t}{2}&0&0&0&c+\frac{t}{2}&0&0&0&0\\
        0&0&r+\frac{t}{4}&0&\frac{t}{4}&0&\frac{t}{4}&c+\frac{t}{4}&0&0\\
        0&0&0&r+\frac{t}{4}&0&\frac{t}{2}&0&0&c+\frac{t}{4}&0\\
        0&0&0&0&0&r+\frac{t}{2}&0&0&c+\frac{t}{2}&0\\
        0&0&0&0&0&0&r+\frac{t}{4}& \frac{t}{2}&0&c+\frac{t}{4}\\
        0&0&0&0&0&0&0&0&0&1\\
\end{pmatrix}\]
\[ T =  \begin{pmatrix}
        r+\frac{t}{2}&c+\frac{t}{2}&0&0&0&0&0&0&0\\
        r+\frac{t}{4}&\frac{t}{4}&\frac{t}{4}&0&c+\frac{t}{4}&0&0&0&0\\
        0&r+\frac{t}{4}&\frac{t}{4}&\frac{t}{4}&0&c+\frac{t}{4}&0&0&0\\
        0&0&r+\frac{t}{2}&\frac{t}{4}&0&0&c+\frac{t}{4}&0&0\\
        0&r+\frac{t}{2}&0&0&0&c+\frac{t}{2}&0&0&0\\
        0&0&r+\frac{t}{4}&0&\frac{t}{4}&0&\frac{t}{4}&c+\frac{t}{4}&0\\
        0&0&0&r+\frac{t}{4}&0&\frac{t}{2}&0&0&c+\frac{t}{4}\\
        0&0&0&0&0&r+\frac{t}{2}&0&0&c+\frac{t}{2}\\
        0&0&0&0&0&0&r+\frac{t}{4}& \frac{t}{2}&0\\
\end{pmatrix}\]

\subsection{Numerical Examples}
In this section we model the game on a $7 \times 7$ Toroidal grid with the following proportion of moves $r=0.4$, $c=0.3$, $t=0.3$.  The table below shows the probability of the robber evading the cop through the first 50 rounds of the game as well as the expected game length based on each of the players' possible starting configurations.  The CoCalc code for these calculations is available to the interested reader in Appendix \ref{sec:code_toroidal} 

\begin{center}
    \begin{tabular}{|c|c|c|c|c|c|c|c|c|c|}
    \hline
        Measure & \multicolumn{9}{c|}{Initial State $d$}\\\hline
                &(3,3)&(3,2)&(3,1)&(3,0)&(2,2)&(2,1)&(2,0)&(1,1)&(1,0)\\\hline
        $\G_{50}(d)$ & $0.5480$&$0.5307$&$0.4942$&$0.4520$&$0.4958$&$0.4376$&$0.3731$&$0.3506$&$0.2331$\\\hline
        $E(d)$ & $78.18$&$95.95$&$71.21$&$65.62$&$71.42$&$63.66$&$54.75$&$51.65$&$34.75$\\\hline
    \end{tabular}
\end{center}

\section{Modeling the game where players sober up over time}\label{sec:Sober}
Harris et al.~asked how to model a tipsy cop and drunk robber game when the cop is sobering up over time \cite[Question~6.1]{HIPSS20}. In this section, we model the game where both players begin completely drunk and sober up as time passes. We define our probability of a tipsy move by either player $t$ as a function of $m$ rounds that have passed $t = f(m)$. Additionally, we set the proportion of sober cop moves to be $c=\frac{a}{b}(1-t)$ and the proportion of sober robber moves to be $r=\frac{b-a}{b}(1-t)$, where $a$ and $b$ are any desired integers that determine the proportion of sober moves that is assigned to the cop and to the robber. As we assume the players are sobering up as time passes, we choose to use a function $f$ with the following properties:
\begin{align*}
f(1) &= 1\\
\lim_{m \to \infty} f(m) &= 0
\end{align*}

With these assumptions, the probability of surviving $\M$ rounds, given an initial state $d$, is given by:
\begin{equation} \G_{\M}(d) = e_d \cdot  \left ( \prod_{m=1}^{\M}  T_m \right ) \mathbbm{1}\label{eq:time_sur}\end{equation}
\\
The expected game time is given by the series:
\begin{equation} E(d) = e_d \cdot \left ( \sum_{n=1}^{\N = \infty} \left (\prod_{m=1}^{n-1} T_m \right) \right ) \cdot \mathbbm{1} \label{eq:time_sober}\end{equation}

We conclude this section with some sample calculations in tables. The code for the calculations is attached in Appendix \ref{sec:code_cycle_sober_time}, so the interested reader can change these calculations to model the game for any values of $a,b,t=f(m), \text{ or } \N$ they choose. 

\subsection{Numerical Examples} 
Assuming the game is played on a cycle graph of six nodes Figure \ref{fig:C_6}, its accompanying Markov Chain is given in Figure \ref{fig:Markov_C_6}.

\noindent\textbf{Example 9.1.1 } Let the tipsiness function be given by the formula $\ds t= f(m) = \frac{4}{m+3}$ as it is an example of a function that follows the rules stated above. Then, the following table shows the probability of the game lasting $5$ rounds from varying starting positions and proportions of sober moves allotted to each player.  It also shows estimates for expected game times from each starting position where we use
$$ E(d) =  e_d \cdot \left ( \sum_{n=1}^{\N} \left (\prod_{m=1}^{n-1} T_m \right) \right ) \mathbbm{1}$$
to approximate the infinite sum in Equation $\eqref{eq:time_sober}$. We chose either $\N=1000$ or $\N=3000$ to approximate these values because these values of $\N$ are sufficiently large to approximate the expected game time accurately to at least five digits as we vary the proportion of sober moves allocated to the robber .

\begin{center}
\begin{tabular}{|c|c|c|c|c|c|c|c|c|c|c|c|}
\hline
   Measure & \multicolumn{11}{c|}{Percentage of Sober Moves that are Robber Moves} \\\hline
    $t=\frac{4}{m+3}$&$0\%$ & $10\%$ & $20\%$ & $30\%$ & $40\%$ & $50\%$ & $60\%$ & $70\%$ & $80\%$ & $90\%$ & $100\%$\\\hline
    $\G_{5}(1)$ & $0.1393$ & $0.1709$ & $0.2023$ & $0.2331$ & $0.2630$ & $0.2917$ & $0.3189$ & $0.3444$ & $0.3681$ & $0.3899$ & $0.4095$\\\hline
    $\G_{5}(2)$ & $0.3667$ & $0.4078$ & $0.4484$ & $0.4882$ & $0.5270$ & $0.5646$ & $0.6008$ & $0.6354$ & $0.6684$ & $0.6995$ & $0.7286$\\\hline
    $\G_{5}(3)$ & $0.2786$ & $0.3419$ & $0.4047$ & $0.4663$  & $0.5260$ & $0.5833$ & $0.6378$ & $0.6889$ & $0.7363$ & $0.7798$ & $0.8190$\\\hline
    $E(1)$ & $3.029^*$ & $3.245^*$ & $3.530^*$ & $3.922^*$ & $4.504^*$ & $5.456^*$ & $7.281^*$ & $11.83^*$ & $30.39^*$& $242.8^{**}$ & $\infty$\\\hline
    $E(2)$ & $4.586^*$ & $4.962^*$ & $5.456^*$ & $6.140^*$ & $7.153^*$ & $8.812^*$ & $11.99^*$ & $19.92^*$ & $52.27^*$ & $422.9^{**}$& $\infty$\\\hline
    $E(3)$ & $5.059^*$ & $5.491^*$ & $6.059^*$ & $6.845^*$ & $8.008^*$ & $9.913^*$ & $13.56^*$ & $22.67^*$ & $59.78^*$ & $484.7^{**}$& $\infty$\\\hline
\end{tabular}\\
*Sum calculated to $\N=1000$, **$\N=3000$
\end{center}
The expected game time for $80\%$ was also calculated at $\N=3000$, and found to be the same for at least the first five digits. It seems only at $90\%$ and above is $\N=1000$ not sufficiently large.

\textbf{Example 9.1.2 } Let tipsiness be the function given by the formula $\ds t = f(m) = \frac{4}{2^m+2}$, because it changes exponentially and follows the rules state above. Then, the following table shows the probability of the game lasting $5$ rounds from varying starting positions and proportions of sober moves allocated to each player. It also shows estimates for expected game times from each starting position where we use 
$$ E(d) = e_d \cdot \left ( \sum_{n=1}^{\N} \left (\prod_{m=1}^{n-1} T_m \right) \right ) \mathbbm{1}$$
to approximate the infinite sum in Equation $\eqref{eq:time_sober}$. We chose either $\N=500$ or $\N=1000$ to approximate these values because these values of $\N$ are sufficiently large for their respective proportion of sober moves allocated to the robber to accurately approximate the expected game time to at least five digits.

\begin{center}
\begin{tabular}{|c|c|c|c|c|c|c|c|c|c|c|c|}
\hline
   Measure & \multicolumn{11}{c|}{Percentage of Sober Moves that are Robber Moves} \\\hline
    $t=\frac{4}{2^m+2}$&$0\%$ & $10\%$ & $20\%$ & $30\%$ & $40\%$ & $50\%$ & $60\%$ & $70\%$ & $80\%$ & $90\%$ & $100\%$\\\hline
    $\G_{5}(1)$ & $0.0381$ & $0.0895$ & $0.1435$ & $0.1978$ & $0.2505$ & $0.3000$ & $0.3449$ & $0.384$ & $0.4177$ & $0.4445$ & $0.4649$\\\hline
    $\G_{5}(2)$ & $0.1926$ & $0.2723$ & $0.3511$ & $0.4275$ & $0.5004$ & $0.5688$ & $0.6317$ & $0.6885$ & $0.7388$ & $0.7820$ & $0.8181$\\\hline
    $\G_{5}(3)$ & $0.0763$ & $0.1790$ & $0.2869$ & $0.3956$  & $0.5010$ & $0.6000$ & $0.6899$ & $0.7688$ & $0.8354$ & $0.8890$ & $0.9298$\\\hline
    $E(1)$ & $2.546^*$ & $2.757^*$ & $3.060^*$ & $3.531^*$ & $4.329^*$ & $5.849^*$ & $9.228^*$ & $18.68^*$ & $57.75^{**}$&$283.7^{**}$& $\infty$\\\hline
    $E(2)$ & $3.753^*$ & $4.119^*$ & $4.647^*$ & $5.467^*$ & $6.858^*$ & $9.507^*$ & $15.40^*$ & $31.90^*$ & $100.2^{**}$&$496.6^{**}$& $\infty$\\\hline
    $E(3)$ & $4.093^*$ & $4.513^*$ & $5.120^*$ & $6.062^*$ & $7.659^*$ & $10.70^*$ & $17.46^*$ & $36.36^*$ & $114.5^{**}$&$566.5^{**}$& $\infty$\\\hline
\end{tabular}\\
*Sum calculated to $\N=500$, **$\N=1000$
\end{center}
The expected game time for $70\%$ was also calculated at $\N=1000$ and found to be the same for at least five digits. However, we are certain that for the expected game time at $90\%$ $\N=1000$ is not sufficiently large, yet trying to calculate the sum with any larger $\N$ causes the CoCalc server to timeout.

\section{Modeling the game where tipsiness is a function of distance}\label{sec:Distance}
Mentioning that it may be more biologically realistic, Harris et al. also asked how to model a tipsy cop and drunk robber game where the players' tipsiness is determined by the distance between them \cite[Question~6.6]{HIPSS20}.
In this section, we model the tipsy cop and robber game where the players sober up as they get closer to each other. This models the scenario where the cop's ability to track the robber improves the closer they get, and the robber senses that the cop is on his trail so he moves more deliberately. The transition matrix $T_\delta$ represents the stochastic matrix, where both the cop and robber sober up as they get closer to each other and get more tipsy as the distance between the two increases.  If $\R$ is the maximum distance the two can be apart on a finite graph (or the distance at which the cop calls off the hunt), then we assume that the function $\delta(d)$ has the following properties:
\begin{align*}
\delta(1) &= 0\\
\lim_{d \to \R} \delta(d) &= 1
\end{align*}
Additionally, we set the proportion of sober cop moves to be $ c=\frac{a}{b}(1-t)$ and the proportion of sober robber moves to be $ r=\frac{b-a}{b}(1-t)$, where $a$ and $b$ are any desired integers that determine what fraction of all sober moves is assigned to the cop and to the robber.

\subsection{Linear Increase of Tipsiness}
In this scenario, we choose to use the function $ \delta (d) = \frac{d-1}{\R}$, where $d$ is the distance between the cop and the robber when they start the chase and $\R$ is the maximum distance they can be apart; this value $\R$ is specified to be either the radius of the graph on finite graphs, or the maximum specified distance before the cop calls off the chase on an infinite graph.
The probability of the game lasting at least $\M$ rounds when starting in state $d$ is
$$\G_M(d) = e_d \cdot  T_\delta^{\M} \cdot \mathbbm{1} .$$
The expected game time is 
$$E(d) = e_d  \cdot \left(I - T_\delta \right)^{-1} \cdot \mathbbm{1} .$$

\subsection{Exponential Increase of Tipsiness}
Now, we choose to use the function $\ds \delta (d) = \frac{1 - 1.2^{(1-d)}}{1 + 1.2^{(1-d)}}$, where $d$ is the distance between the cop and the robber when they start the chase. 
The probability of the game lasting through $\M$ rounds is given by:
$$\G_M(d) = e_d \cdot  T_\delta^{\M} \cdot \mathbbm{1} $$

The expected game time is given by:
$$E(d) = e_d  \cdot \left(I - T_\delta \right)^{-1} \cdot \mathbbm{1} $$

We conclude this section with some sample calculations in tables. We have published the accompanying CoCalc code for these computations in Appendix \ref{sec:code_tree_exponential_tipsiness}, so the interested reader can change these calculations to model the game for any values of $\n, m, d, r$ that they choose. 

\subsection{Numerical Example on Cycle Graph}
Given a cycle graph of $\n=10$ nodes, we will have the possibilities of starting at distances $0$ to $5$ away with maximum distance $\R = \n/2 = 5$. The accompanying Markov Chain is (again we assume that the cop always moves and the sober robber chooses not to move if they are furthest away from each other as noted in Section \ref{sec:Cycle Graphs} - even nodes cycle graphs) 
\begin{center}
 \begin{tikzpicture}
        \node[state]             (0) {0};
        \node[state, right=of 0] (1) {1};
        \node[state, right=of 1] (2) {2};
        \node[state, right=of 2] (3) {3};
        \node[state, right=of 3] (4) {4};
        \node[state, right=of 4] (5) {5};

        \draw[every loop]
            (1) edge[bend right, auto=right] node {$c_1 + \frac{t_1}{2}$} (0)
            (0) edge[loop above, auto=right] node {1} (1)
            
            (2) edge[bend right, auto=right] node {$c_2 + \frac{t_2}{2}$} (1)
            (1) edge[bend right, auto=right] node {$r_1 + \frac{t_1}{2}$} (2)
            (2) edge[bend right, auto=right] node {$r_2 + \frac{t_2}{2}$} (3)
            (3) edge[bend right, auto=right] node {$c_3 + \frac{t_3}{2}$} (2)
            (3) edge[bend right, auto=right] node {$r_3 + \frac{t_3}{2}$} (4)
            (4) edge[bend right, auto=right] node {$c_4 + \frac{t_4}{2}$} (3)
            (4) edge[bend right, auto=right] node {$r_4 + \frac{t_4}{2}$} (5)
            (5) edge[loop above, auto=right] node {$r_5$} (5)
            (5) edge[bend right, auto=right ] node {$c_5 + t_5$} (4);
    \end{tikzpicture}
\end{center}

Based on the Markov chain we get the transformation matrix based on distance $T_\delta$
\[ T_\delta= \begin{pmatrix}
    0 & r_1 + \frac{t_1}{2} & 0 & 0 & 0\\
    c_2 + \frac{t_2}{2} & 0 & r_2 + \frac{t_2}{2} & 0 & 0\\
    0 & c_3 + \frac{t_3}{2} & 0 & r_3 + \frac{t_3}{2} & 0\\
    0 & 0 & c_4+ \frac{t_4}{2} & 0 & r_4 + \frac{t_4}{2}\\
    0 & 0 & 0 & c_5 + t_5 & r_5 \\
\end{pmatrix}
\]
Note that for each distance $d$, $r_d+c_d+t_d=1$. In our model we specify what percentage of all sober moves $(1-t_d)$ are given to the robber.

The probability of surviving 20 rounds when starting at distance $d$ apart $\G_{20}(d)$ and the expected duration of the chase $E(d)$ using the linear growth of the tipsiness $\ds t = \delta (d) = \frac{d-1}{\R}$ is shown in the table below.

\begin{center}
\begin{tabular}{|c|c|c|c|c|c|c|c|c|c|c|c|}
\hline
   Measure & \multicolumn{11}{c|}{Percentage of Sober Moves that are Robber Moves} \\\hline
     & $0\%$ & $10\%$ & $20\%$ & $30\%$ & $40\%$ & $50\%$ & $60\%$ & $70\%$ & $80\%$ & $90\%$ & $100\%$\\\hline
    $\G_{20}(1)$ & $0$ & $0.0001$ & $0.0024$ & $0.015$ & $0.057$ & $0.148$ & $0.299$ & $0.491$ & $0.689$ & $0.864$ & $1$\\\hline
    $\G_{20}(2)$ & $6E-5$ & $0.001$ & $0.009$ & $0.041$ & $0.125$ & $0.277$ & $0.48$ & $0.69$ & $0.86$ & $0.958$ & $1$\\\hline
    $\G_{20}(3)$ & $0.0005$ & $0.0049$ & $0.025$ & $0.085$  & $0.208$ & $0.391$ & $0.599$ & $0.784$ & $0.912$ & $0.978$ & $1$\\\hline
    $\G_{20}(4)$ & $0.0011$ & $0.0088$ & $0.039$ & $0.118$ & $0.261$ & $0.454$ & $0.657$ & $0.824$ & $0.932$ & $0.984$ & $1$\\\hline
    $\G_{20}(5)$ & $0.0024$ & $0.0139$ & $0.053$ & $0.144$ & $0.295$ & $0.489$ & $0.683$ & $0.839$ & $0.938$ & $0.985$ & $1$\\\hline
    $E(1)$ & $1$ & $1.29$ & $1.81$ & $2.78$ & $4.77$ & $9.25$ & $20.58$ & $53.89$ & $177.7$ & $914.8$ & $\infty$\\\hline
    $E(2)$ & $2.30$ & $2.97$ & $4.06$ & $5.93$ & $9.41$ & $16.5$ & $32.64$ & $75.55$ & $220.8$ & $1015$ & $\infty$\\\hline
    $E(3)$ & $4.02$ & $5.04$ & $6.59$ & $9.1$ & $13.45$ & $21.75$ & $39.64$ & $85.2$ & $234.7$ & $1036$ & $\infty$\\\hline
    $E(4)$ & $5.87$ & $7.09$ & $8.88$ & $11.65$ & $16.32$ & $25.0$ & $43.36$ & $89.49$ & $239.7$ & $1042$ & $\infty$\\\hline
    $E(5)$ & $6.87$ & $8.14$ & $9.97$ & $12.79$ & $17.51$ & $26.25$ & $44.68$ & $90.88$ & $241.2$ & $1044$ & $\infty$\\\hline
\end{tabular}\\
\end{center}

Similarly, the probability of surviving 20 rounds when starting distance $d$ apart $\G_{20}(d)$ and the expected duration of the chase $E(d)$ using the exponential growth for tipsiness $\ds t = \delta (d) = \frac{1 - 1.2^{(1-d)}}{1 + 1.2^{(1-d)}}$ is shown in the table below.

\begin{center}
\begin{tabular}{|c|c|c|c|c|c|c|c|c|c|c|c|}
\hline
   Measure & \multicolumn{11}{c|}{Percentage of Sober Moves that are Robber Moves} \\\hline
    &$0\%$ & $10\%$ & $20\%$ & $30\%$ & $40\%$ & $50\%$ & $60\%$ & $70\%$ & $80\%$ & $90\%$ & $100\%$\\\hline
    $\G_{20}(1)$ & $0$ & $1.2E-5$ & $0.0008$ & $0.0093$ & $0.049$ & $0.152$ & $0.325$ & $0.532$ & $0.726$ & $0.881$ & $1$\\\hline
    $\G_{20}(2)$ & $4E-8$ & $6.1E-5$ & $0.0025$ & $0.025$ & $0.109$ & $0.287$ & $0.529$ & $0.755$ & $0.906$ & $0.978$ & $1$\\\hline
    $\G_{20}(3)$ & $9E-7$ & $0.0003$ & $0.0072$ & $0.0513$  & $0.181$ & $0.404$ & $0.656$ & $0.852$ & $0.957$ & $0.993$ & $1$\\\hline
    $\G_{20}(4)$ & $1.9E-6$ & $0.0006$ & $0.012$ & $0.0741$ & $0.234$ & $0.477$ & $0.722$ & $0.891$ & $0.972$ & $0.997$ & $1$\\\hline
    $\G_{20}(5)$ & $8.2E-6$ & $0.0012$ & $0.018$ & $0.0947$ & $0.269$ & $0.515$ & $0.749$ & $0.904$ & $0.976$ & $0.997$ & $1$\\\hline
    $E(1)$ & $1$ & $1.27$ & $1.72$ & $2.59$ & $4.52$ & $9.58$ & $25.7$ & $90.6$ & $466.9$ & $4825$ & $\infty$\\\hline
    $E(2)$ & $2.1$ & $2.67$ & $3.59$ & $5.29$ & $8.79$ & $17.16$ & $41.1$ & $128$ & $582.3$ & $5360$ & $\infty$\\\hline
    $E(3)$ & $3.33$ & $4.19$ & $5.56$ & $7.93$ & $12.52$ & $22.74$ & $50.1$ & $144$ & $615$ & $5443$ & $\infty$\\\hline
    $E(4)$ & $4.64$ & $5.74$ & $7.41$ & $10.18$ & $15.31$ & $26.31$ & $54.8$ & $150.6$ & $624.8$ & $5459$ & $\infty$\\\hline
    $E(5)$ & $5.64$ & $6.82$ & $8.58$ & $11.46$ & $16.73$ & $27.89$ & $56.6$ & $152.7$ & $627.2$ & $5462$ & $\infty$\\\hline
\end{tabular}\\
\end{center}

Based on the results in the tables the chase is longer if the tipsiness is increasing exponentially and the percentage of all sober moves $(1-t_d)$ are given to the robber is less than $50\%$. In the case that the percentage of all sober moves $(1-t_d)$ are given to the robber is greater than $50\%$ the chase is longer if the tipsiness is increasing linearly. Assuming the robber is given the choice, the robber should pick a strategy based on the percentage of sober moves assigned to him--if he gets to move more than $50\%$ of the sober moves his tipsiness should increase linearly, otherwise (he gets to move $50\%$ or less of all sober moves) he has a better chance of surviving if the tipsiness changes exponentially.

\subsection{Infinite Regular Trees}

If we play the game on an infinite regular tree, and the cop calls off the hunt if the robber reaches a distance of $\R$ nodes away from the cop, then the matrix $P$ is an $(\R+1 \times \R+1)$ tridiagonal matrix with $P_{1,1} = 1$, $P_{\R+1,\R+1} = 1$ and the upper and lower diagonal entries are $P_{k+1,k} = c_{k+1} + \frac{t_{k+1}}{\Delta}$ for $1 \leq k \leq \R-1$, and    $P_{k+1,k+2} = t_{k+1} \frac{\Delta - 1}{\Delta} +r_{k+1}$ for $1 \leq k \leq \R-1$, and finally $P_{i,j}=0$ otherwise.  To obtain the matrix $T$ we remove from $P$ the first and last rows and columns corresponding to both absorbing states. 

\begin{figure}[h]
\begin{center}
    \begin{tikzpicture}[scale=0.75]
        \node[state] at (0,0) (0) {$0$};
        \node[state] at (3,0) (1) {$1$};
        \node[state] at (6,0) (2) {$2$};
        \node[state] at (9,0) (3) {$3$};
        \node[state] at (12,0) (4) {$4$};
        \node[state] at (15,0) (5) {5};
        \draw[every loop]
            (1) edge[bend right, auto=right] node {$c_1 + \frac{t_1}{\Delta}$} (0)
            (0) edge[loop above] node {$1$} (0)
            (5) edge[loop above] node {$1$} (5)
            (1) edge[bend right, auto=right] node {$r_1 + t_1\frac{\Delta - 1}{\Delta}$} (2)
            (2) edge[bend right, auto=right] node {$r_2 + t_2\frac{\Delta - 1}{\Delta}$} (3)
            (3) edge[bend right, auto=right] node {$r_3 + t_3\frac{\Delta - 1}{\Delta}$} (4)
            (4) edge[bend right, auto=right] node {$r_4 + t_4\frac{\Delta - 1}{\Delta}$} (5)
            (2) edge[bend right, auto=right] node {$c_2 + \frac{t_2}{\Delta}$} (1)
            (3) edge[bend right, auto=right] node {$c_3 + \frac{t_3}{\Delta}$} (2)
            (4) edge[bend right, auto=right] node {$c_4 + \frac{t_4}{\Delta}$} (3);
    \end{tikzpicture}
\end{center}
\captionof{figure}{Markov chain for regular tree with maximum distance of $\R=5$.} \label{fig:rt5}
\end{figure}
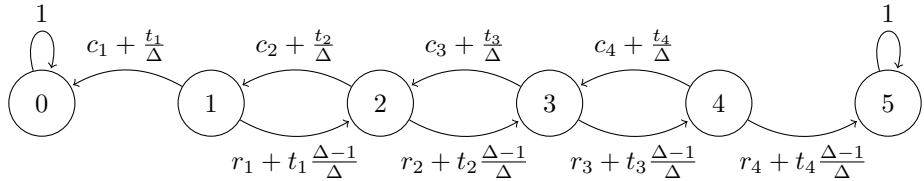

For example, on an infinite regular tree of degree $\Delta$ with maximum specified distance $\R =5$, the accompanying Markov Chain is in Figure \ref{fig:rt5}.
Based on the Markov chain we get the transformation matrix based on distance $T_\delta$
 $$T_\delta= \begin{pmatrix}
    0 & r_1 + t_1\frac{\Delta - 1}{\Delta} & 0 & 0 \\
    c_2 + \frac{t_2}{\Delta} & 0 & r_2 + t_2\frac{\Delta - 1}{\Delta} & 0\\
    0 & c_3 + \frac{t_3}{\Delta} & 0 & r_3 + t_3\frac{\Delta - 1}{\Delta}\\
    0 & 0 & c_4 + \frac{t_4}{\Delta} & 0 \\
\end{pmatrix}.$$
Note that in this case we remove states $0$ and $5$ because they are absorbing.  Similarly in the numerical example on cycle graphs, we specify what percentage of all sober moves $(1-t_d)$ are given to the robber.

\subsection{Numerical Example on Infinite Regular Tree with $\Delta = 4$ and Maximum Distance $\R = 10$}
The probability that the game continues for $\M= 30$ rounds when starting distance $d$ apart $\G_{30}(d)$ and the expected duration of the chase $E(d)$ when $\Delta = 4$ when tipsiness increases linearly $\ds t = \delta (d) = \frac{d-1}{\R}$ is shown in the table below.  Note that we assume the cop calls off the chase if the distance reaches $\R=10$.

\begin{center}
\begin{tabular}{|c|c|c|c|c|c|c|c|c|c|c|c|}
\hline
   Measure & \multicolumn{11}{c|}{Percentage of Sober Moves that are Robber Moves} \\\hline
    & $0\%$ & $10\%$ & $20\%$ & $30\%$ & $40\%$ & $50\%$ & $60\%$ & $70\%$ & $80\%$ & $90\%$ & $100\%$\\\hline
    $\G_{30}(1)$ & $0$ & $4.5E-5$ & $0.0008$ & $0.006$ & $0.022$ & $0.045$ & $0.051$ & $0.033$ & $0.012$ & $0.002$ & $0.0002$\\\hline
    $\G_{30}(2)$ & $1.4E-5$ & $0.0003$ & $0.0031$ & $0.016$ & $0.046$ & $0.073$ & $0.067$ & $0.034$ & $0.009$ & $0.0013$ & $0.0001$ \\\hline
    $\G_{30}(3)$ & $0.0002$ & $0.0017$ & $0.0099$ & $0.035$  & $0.077$ & $0.099$ & $0.077$ & $0.034$ & $0.008$ & $0.0012$ & $9.7E-5$ \\\hline
    $\G_{30}(4)$ & $0.0007$ & $0.0042$ & $0.0174$ & $0.048$ & $0.083$ & $0.088$ & $0.057$ & $0.021$ & $0.004$ & $0.0005$ & $3.9E-5$\\\hline
    $\G_{30}(5)$ & $0.0026$ & $0.0102$ & $0.0303$ & $0.063$ & $0.089$ & $0.08$ & $0.046$ & $0.016$ & $0.003$ & $0.0004$ & $3.5E-5$\\\hline
    $\G_{30}(6)$ & $0.0041$ & $0.0128$ & $0.031$ & $0.054$ & $0.065$ & $0.051$ & $0.025$ & $0.008$ & $0.001$ & $0.0002$ & $1.3E-5$\\\hline
    $\G_{30}(7)$ & $0.0065$ & $0.0161$ & $0.0318$ & $0.046$ & $0.047$ & $0.033$ & $0.015$ & $0.005$ & $0.0009$ & $0.0001$ & $9.4E-6$\\\hline
    $\G_{30}(8)$ & $0.0045$ & $0.0102$ & $0.0182$ & $0.024$ & $0.022$ & $0.014$ & $0.006$ & $0.0016$ & $0.0003$ & $3.3E-5$ & $2.6E-6$\\\hline
    $\G_{30}(9)$ & $0.0027$ & $0.0054$ & $0.0086$ & $0.01$ & $0.009$ & $0.005$ & $0.002$ & $0.0006$ & $0.0001$ & $1.4E-5$ & $1.2E-6$\\\hline
    $E(1)$ & $1$ & $1.29$ & $1.82$ & $2.86$ & $4.75$ & $7.3$ & $9.58$ & $11$ & $11.5$ & $11.6$ & $11.4$\\\hline
    $E(2)$ & $2.2$ & $2.88$ & $4.09$ & $6.21$ & $9.38$ & $12.6$ & $14.3$ & $14.2$ & $13.2$ & $11.8$ & $10.4$\\\hline
    $E(3)$ & $3.72$ & $4.9$ & $6.8$ & $9.68$ & $13.1$ & $15.5$ & $15.6$ & $14.2$ & $12.3$ & $10.7$ & $9.39$\\\hline
    $E(4)$ & $5.64$ & $7.29$ & $9.62$ & $12.5$ & $15.15$ & $16$ & $14.8$ & $12.8$ & $10.9$ & $9.38$ & $8.28$\\\hline
    $E(5)$ & $7.81$ & $9.61$ & $11.77$ & $13.9$ & $15.19$ & $14.7$ & $12.8$ & $10.8$ & $9.17$ & $7.99$ & $7.11$\\\hline
    $E(6)$ & $9.54$ & $10.97$ & $12.37$ & $13.4$ & $13.4$ & $12.1$ & $10.3$ & $8.63$ & $7.41$ & $6.53$ & $5.87$\\\hline
    $E(7)$ & $9.75$ & $10.44$ & $10.92$ & $10.9$ & $10.3$ & $8.99$ & $7.58$ & $6.44$ & $5.61$ & $5.01$ & $4.55$\\\hline
    $E(8)$ & $7.8$ & $7.86$ & $7.75$ & $7.37$ & $6.67$ & $5.76$ & $4.91$ & $4.25$ & $3.78$ & $3.42$ & $3.14$\\\hline
    $E(9)$ & $4.11$ & $3.99$ & $3.79$ & $3.51$ & $3.14$ & $2.73$ & $2.37$ & $2.11$ & $1.91$ & $1.75$ & $1.63$\\\hline
\end{tabular}\\
\end{center}

The table below shows the probability of a game lasting 30 rounds when starting at distance $d$ apart $\G_{30}(d)$ and the expected duration of the chase $E(d)$ when the game takes place on a regular tree of degree $\Delta =4$ and tipsiness is modeled using the exponential function $\ds t_d = \delta (d) = \frac{1 - 1.2^{1-d}}{1 + 1.2^{1-d}}$.
Again we assume the cop calls off the chase if the distance reaches $\R=10$.

\begin{center}
\begin{tabular}{|c|c|c|c|c|c|c|c|c|c|c|c|}
\hline
   Measure & \multicolumn{11}{c|}{Percentage of Sober Moves that are Robber Moves} \\\hline
   & $0\%$ & $10\%$ & $20\%$ & $30\%$ & $40\%$ & $50\%$ & $60\%$ & $70\%$ & $80\%$ & $90\%$ & $100\%$\\\hline
   $\G_{30}(1)$ & $0$ & $0.0004$ & $0.0016$ & $0.0042$ & $0.0085$ & $0.014$ & $0.019$ & $0.021$ & $0.02$ & $0.017$ & $0.014$\\\hline
    $\G_{30}(2)$ & $0.0014$ & $0.0029$ & $0.0056$ & $0.0098$ & $0.015$ & $0.02$ & $0.022$ & $0.021$ & $0.017$ & $0.013$ & $0.0089$ \\\hline
    $\G_{30}(3)$ & $0.0055$ & $0.0082$ & $0.012$ & $0.017$  & $0.021$ & $0.024$ & $0.024$ & $0.021$ & $0.016$ & $0.012$ & $0.0085$ \\\hline
    $\G_{30}(4)$ & $0.0066$ & $0.0084$ & $0.0107$ & $0.013$ & $0.015$ & $0.015$ & $0.014$ & $0.012$ & $0.009$ & $0.007$ & $0.0048$\\\hline
    $\G_{30}(5)$ & $0.0073$ & $0.0086$ & $0.0099$ & $0.011$ & $0.012$ & $0.012$ & $0.01$ & $0.0087$ & $0.007$ & $0.005$ & $0.0039$\\\hline
    $\G_{30}(6)$ & $0.0041$ & $0.0046$ & $0.005$ & $0.0054$ & $0.006$ & $0.005$ & $0.005$ & $0.0039$ & $0.003$ & $0.002$ & $0.0018$\\\hline
    $\G_{30}(7)$ & $0.0027$ & $0.003$ & $0.003$ & $0.0033$ & $0.003$ & $0.003$ & $0.0027$ & $0.0023$ & $0.0018$ & $0.0015$ & $0.0012$\\\hline
    $\G_{30}(8)$ & $0.001$ & $0.001$ & $0.001$ & $0.0012$ & $0.0011$ & $0.001$ & $0.0009$ & $0.0008$ & $0.0007$ & $0.0005$ & $0.0004$\\\hline
    $\G_{30}(9)$ & $0.004$ & $0.0004$ & $0.0004$ & $0.0004$ & $0.0004$ & $0.0004$ & $0.0003$ & $0.0003$ & $0.0002$ & $0.0002$ & $0.0002$\\\hline
    $E(1)$ & $1$ & $1.47$ & $2.22$ & $3.32$ & $4.8$ & $6.6$ & $8.55$ & $10.5$ & $12.1$ & $13.5$ & $14.5$\\\hline
    $E(2)$ & $3.6$ & $4.7$ & $6.09$ & $7.73$ & $9.5$ & $11.2$ & $12.6$ & $13.5$ & $13.9$ & $13.9$ & $13.5$\\\hline
    $E(3)$ & $7.41$ & $8.53$ & $9.71$ & $10.9$ & $12$ & $12.8$ & $13.2$ & $13.3$ & $13.1$ & $12.8$ & $12.4$\\\hline
    $E(4)$ & $9.85$ & $10.5$ & $11$ & $11.5$ & $11.8$ & $12.1$ & $12.1$ & $11.9$ & $11.6$ & $11.3$ & $11$\\\hline
    $E(5)$ & $9.88$ & $10.1$ & $10.2$ & $10.3$ & $10.4$ & $10.3$ & $10.2$ & $10$ & $9.81$ & $9.6$ & $9.4$\\\hline
    $E(6)$ & $8.38$ & $8.4$ & $8.4$ & $8.38$ & $8.33$ & $8.26$ & $8.16$ & $8.04$ & $7.91$ & $7.79$ & $7.67$\\\hline
    $E(7)$ & $6.33$ & $6.31$ & $6.28$ & $6.25$ & $6.21$ & $6.15$ & $6.09$ & $6.02$ & $5.96$ & $5.89$ & $5.83$\\\hline
    $E(8)$ & $4.18$ & $4.17$ & $4.15$ & $4.13$ & $4.1$ & $4.07$ & $4.04$ & $4.01$ & $3.98$ & $3.95$ & $3.92$\\\hline
    $E(9)$ & $2.07$ & $2.06$ & $2.05$ & $2.05$ & $2.04$ & $2.03$ & $2.02$ & $2$ & $1.99$ & $1.98$ & $1.97$\\\hline
\end{tabular}\\
\end{center}

Based on the results in the tables above, the robber should pick the strategy where the tipsiness changes exponentially, since it gives him a higher chance of survival and the chase lasts longer.

\section{Future Directions and Open Questions}

We end this paper by pointing out a few possible areas of future study and posing a few open questions:
Modeling the game on infinite grids is considerably more difficult than on infinite regular trees because there are multiple states for a given distance, and it is not immediately clear what strategy the cop should use to decrease distance or the robber should use to increase distance.

\begin{enumerate}
\item On an infinite grid, what are the optimal strategies for the cop and robber?  
\item Who is more likely to win on an infinite grid, if both players employ their optimal strategies?
\item How will the game change if the cop and the robber do not take turns, but instead move simultaneously?   
\end{enumerate}

\section{Acknowledgements}

We would like to thank Drs. Pamela Harris, Florian Lehner, and Alicia Prieto-Langarica for many helpful and inspiring conversations and suggestions regarding this project. We also thank Drs. Katie Johnson and Shaun Sullivan for reading and providing helpful feedback regarding earlier drafts of this manuscript. This collaboration was supported by the FGCU Seidler Collaboration Fellowship.

\bibliography{Bibliography}
\bibliographystyle{plain}

\newpage
\appendix
\section{SageMath code}\label{sec:code}
The code described in the previous Sections is provided in this appendix.

\subsection{Sage Code for Petersen Graph}\label{sec:code_Petersen}
\begin{lstlisting}
var ('c,t,M,i,p')

def check_rt(r,t):
    check=1
    if r<0:
        check=0
        print("r cannot be negative")
    if t<0:
        check=0
        print("t cannot be negative")
    if t+r>1:
        check=0
        print("r+t cannot be greater than 1")
    return check

def define_P(n=4): #default n=4 can be changed later
    L=[]
    for i in range(n+1):
        L.append([])
        for j in range(n+1):
            L[i].append(0)
    L[0][0]=1
    L[-1][-1]=1
    for i in range(1,n):
        L[i][i+1]=r+(2*t/3)
    for i in range(1,n):
        L[i][i-1] =c+t/3
    L[n][n]=r+(2*t/3)
    L[n][n-1]=c+t/3
    M=matrix(L)
    return M
######################################

def check_define_P(n=4):
    if check_rt(r,t) ==1:
        return define_P(n)

M=7 #number of rounds
r=.30
t=.40
c=1-t-r
n=10 #default n=4 can be changed here
check_define_P(n)

######################################
def define_T(n):
    P=define_P(n)
    T=P.submatrix(1,1,n,n)
    return T
######################################

P=define_P(n)
print(P)
print()

T=define_T(n)
print(T)
print()

I = identity_matrix(n)
print(I)


one=[]
for i in range(n):
    one.append(1)
N=vector(one)

print(N)

print("c=",c,"r=",r)
for d in range(1,n+1):
    e_d = I.column(d-1)
    print("d=",d,"n=",n,"M=",M)
    G = e_d*(T^M)*N
    print("Survival probability of M rounds from distance ", d, " is ", G)
    E = e_d*((I-T)^(-1))*N
    print("Expected number of rounds from distance ", d, " is " , E)
\end{lstlisting}

\newpage
\subsection{Sage Code for $7 \times 7$ Toroidal Grid}\label{sec:code_toroidal}
\begin{lstlisting}
var('r,t,c')
r=0.4 #Proportion of moves that are sober robber moves.
c=0.3 #Proportion of moves that are sober cop moves.
t=0.3 #Proportion of moves that are tipsy moves by either player.
M=50  #Number of rounds you want the robber to survive.
#####################################################
if c+r+t != 1:                                  #
    print 'Make sure your probabilities sum to 1:'  # Probabilities Check
    print 'Sum = ' +`t+r+c` + ' = 1 ?'           #
#####################################################
else:
    v=vector([1,1,1,1,1,1,1,1,1])
    a=vector([1,0,0,0,0,0,0,0,0]) #starting state (3,3)
    b=vector([0,1,0,0,0,0,0,0,0]) #starting state (3,2)
    d=vector([0,0,1,0,0,0,0,0,0]) #starting state (3,1)
    e=vector([0,0,0,1,0,0,0,0,0]) #starting state (3,0)
    f=vector([0,0,0,0,1,0,0,0,0]) #starting state (2,2)
    g=vector([0,0,0,0,0,1,0,0,0]) #starting state (2,1)
    h=vector([0,0,0,0,0,0,1,0,0]) #starting state (2,0)
    i=vector([0,0,0,0,0,0,0,1,0]) #starting state (1,1)
    j=vector([0,0,0,0,0,0,0,0,1]) #starting state (1,0)
    T=matrix([[r+t/2,c+t/2,0,0,0,0,0,0,0],
          [r+t/4,t/4,t/4,0,c+t/4,0,0,0,0],
          [0,r+t/4,t/4,t/4,0,c+t/4,0,0,0],
          [0,0,r+t/2,t/4,0,0,c+t/4,0,0],
          [0,r+t/2,0,0,0,c+t/2,0,0,0],
          [0,0,r+t/4,0,t/4,0,t/4,c+t/4,0],
          [0,0,0,r+t/4,0,t/2,0,0,c+t/4],
          [0,0,0,0,0,r+t/2,0,0,c+t/2],
          [0,0,0,0,0,0,r+t/4,t/2,0]])
    G=vector((T^M)*v)
    E=vector((1-T).inverse()*v)
    print 'G_'+`M`+'(3,3)=',G*a
    print 'G_'+`M`+'(3,2)=',G*b
    print 'G_'+`M`+'(3,1)=',G*d
    print 'G_'+`M`+'(3,0)=',G*e
    print 'G_'+`M`+'(2,2)=',G*f
    print 'G_'+`M`+'(2,1)=',G*g
    print 'G_'+`M`+'(2,0)=',G*h
    print 'G_'+`M`+'(1,1)=',G*i
    print 'G_'+`M`+'(1,0)=',G*j
    print 'E(3,3)=',E*a
    print 'E(3,2)=',E*b
    print 'E(3,1)=',E*d
    print 'E(3,0)=',E*e
    print 'E(2,2)=',E*f
    print 'E(2,1)=',E*g
    print 'E(2,0)=',E*h
    print 'E(1,1)=',E*i
    print 'E(1,0)=',E*j
\end{lstlisting}

\newpage
\subsection{Sage Code for Cycle Graphs Sobering up Overtime}\label{sec:code_cycle_sober_time}
\begin{lstlisting}
from sage.calculus.calculus import symbolic_sum
var ('t,M,N,i,p')
M=5     #number of rounds you want the robber to survive
N=1000  #Nth partial sum of infinite sum. Higher is more accurate but takes longer.
k=9     #number of nodes in cycle graph
p=0.51  #proportion of sober moves allocated to robber
#####################################################
#    if p > 0.8 then N should be at at least 3000.  #
#       otherwise N=1000 is sufficient.             #
#####################################################


if k%2==0:
        n=(k/2)
else: n=(k-1)/2

def define_P(n,p,t):
    L=[]
    for i in range(n+1):
        L.append([])
        for j in range(n+1):
            L[i].append(0)
    L[0][0]=1
    L[-1][-1]=1
    for i in range(1,n):
        L[i][i+1]=p*(1-t) +t/2 # transition probability to make distance bigger
        L[i][i-1]=(1-p)*(1-t) + t/2 # transition probability to make distance smaller
    if k%2==0: #if cycgraph is even
        L[n][n]= p*(1-t)
        L[n][n-1]=(1-p)*(1-t)+t
    else: #if cycle graph is odd
        L[n][n]=p*(1-t)+t/2
        L[n][n-1]=(1-p)*(1-t)+t/2

    M=matrix(L)
    return M
###################################################################################

def define_T(n,p,t):
    P=define_P(n,p,t)
    T=P.submatrix(1,1,n,n)
    return T

####################################################################################
# ouptut: n length vector of all ones
#####################|###############################################################
def one_vector(n):
    L=[]
    for i in range(n):
        L.append(1)
    return vector(L)
####################################################################################
def position_vector(n,d):
    L=[]
    for i in range(n):
        L.append(0)
    L[d]=1
    return vector(L)
####################################################################################


T=define_T(n,p,t)
#print T
#Uncomment the line above to see Transition Matrix

print 'Calculating...'
x=vector(prod(define_T(n,p,4/(3+m)) for m in range (1,M+1))*one_vector(n))
y=vector(sum(prod(define_T(n,p,4/(3+m)) for m in range (1,q))for q in range (1,N+1))*one_vector(n))
print 'done:'

for i in range (0,n):
    print 'G_'+`M`+'('+`i+1`+')=',x*position_vector(n,i)
    print 'E('+`i+1`+')=',y*position_vector(n,i)
\end{lstlisting}

\newpage
\subsection{Sage Code for Friendship Graphs}\label{sec:code_friendship}
\begin{lstlisting}
var('r,tr,tc,n,c')
tr=.4 #Probability of a tipsy robber move
tc=.4 #Probability of a tipsy cop move
r=.1  #Probability of a sober robber move
c=.1  #Probability of a sober cop move
n=6   #Number of triangles
M=5  #Number of rounds the game lasts
#####################################################
if tc+tr+r+c != 1:                                  #
    print 'Make sure your probabilities sum to 1:'  # Probabilities Check
    print 'Sum = '+`tc+tr+r+c` + ' = 1 ?'           #
#####################################################
else:
    T=matrix([[r+tc/2+tr/2,c+tc/2,tr/2,0],
              [tc*((n-1)/n),r+tr/2,0,tc/(2*n)],
              [r+tr*((n-1)/n),0,tc/2,tr/(2*n)],
              [0,tc/2,r+tr/2,0]])
    v=vector([1,1,1,1])
    a=vector([1,0,0,0]) #starting state 2
    b=vector([0,1,0,0]) #starting state 1cc
    c=vector([0,0,1,0]) #starting state 1rc
    d=vector([0,0,0,1]) #starting state 1e
    E=vector((matrix.identity(4)-T).inverse()*v)
####################################################
    print 'G_'+`M`+'(2)=',a*T^M*v
    print 'G_'+`M`+'(1cc)=',b*T^M*v
    print 'G_'+`M`+'(1rc)=',c*T^M*v
    print 'G_'+`M`+'(1e)=',d*T^M*v
    print 'E(2)=',E*a
    print 'E(1cc)=',E*b
    print 'E(1rc)=',E*c
    print 'E(1e)=',E*d
\end{lstlisting}

\newpage
\subsection{Sage Code for Regular Tree with Exponential Change of Tipsiness as a Function of Distance}\label{sec:code_tree_exponential_tipsiness}
\begin{lstlisting}
var ('c,t,M,i,p')
M=30 #number of rounds you want the robber to survive
#c=1-t-c #proportion of sober cop moves
#r--proportion of sober robber moves
#t--proportion of tipsy moves, changes with distance
#D--distance between cop and robber when they start
n=10#maximum distance between cop and robber
delta=4
p=1 #proportion of sober moves allocated to robber


###############################################################################
# input: n is maximum distance between cop and robber before calling off chase
# input: p is proportion of sober moves allocated to robber 
# output: transition matrix (including 2 absorbing states)
###############################################################################
def define_P(n,p):
    L=[]
    for i in range(n+1):
        L.append([])
        for j in range(n+1):
            L[i].append(0)
    T=[0]
    R=[0]
    C=[0]
    L[0][0]=1
    L[-1][-1]=1
    for D in range(1,n+1):
        s=1 - (1.2)^(1-D)
        q=1 + (1.2)^(1-D)
        T.append(s/q) #tipsiness based on distance
        R.append(p*(1-T[D])) #Percentage of sober moves are robber moves as a function of t based on distance
        C.append(1-T[D]-R[D])
    print "T=",T
    print "R=",R
    print "C=",C
    for i in range(1,n):
        L[i][i+1]=R[i]+(T[i]*((delta - 1)/delta)) # transition probability to make distance bigger
        #print "L[",i,"][",i+1,"]=",L[i][i+1]
        L[i][i-1] = C[i]+T[i]*(1/delta) # transition probability to make distance smaller
    L[n][n-1]=0
    L[n][n]=1
    M=matrix(L)
    return M
###################################################################################

####################################################################################
# input: n is maximum distance between cop and robber before calling off chase
# input: p is proportion of sober moves allocated to robber 
# output: transition matrix (without absorbing states)
####################################################################################
def define_T(n,p):
    P=define_P(n,p)
    T=P.submatrix(1,1,n-1,n-1)
    return T
####################################################################################

####################################################################################
# input: n the maximum distance between cop and robber before calling off chase
# ouptut: n-2 length vector of all ones 
####################################################################################
def one_vector(n):
    L=[]
    for i in range(n-1):
        L.append(1)
    return vector(L)
####################################################################################

####################################################################################
#input: n the maximum distance between cop and robber before calling off chase
#input: d such that d+1 is starting distance between cop and robber
####################################################################################
def position_vector(n,d):
    L=[]
    for i in range(n-1):
        L.append(0)
    L[d]=1
    return vector(L)
####################################################################################


P=define_P(n,p)
#print P
#print

T=define_T(n,p)
#print T

#print

one= one_vector(n) 

print 
print 'G_d is probability of surviving', M, 'rounds starting at distance d.'
print 'E_d is expected game length E_d from distance d:'

print 
for d in range(n-1):
    print "starting distance d=", d+1
    e=position_vector(n,d) 
    G_d=e*(T^M)*one
    print 'G_d=', float(G_d)
    E_d=e*((1-T).inverse())*one
    print 'E_d=' ,float(E_d)
    print 
\end{lstlisting}
\end{document}